\documentclass{amsart}

\usepackage[linktoc=page]{hyperref}

\usepackage[utf8]{inputenc}
\usepackage{cite}
\usepackage{amssymb}
\usepackage{amsthm}
\usepackage{amsfonts}
\usepackage{amsmath}
\usepackage{mathrsfs}
\usepackage[all]{xy}
\usepackage{graphicx}
\usepackage{mathrsfs}
\usepackage{extpfeil}
\usepackage{mathtools}
\usepackage{tikz-cd}
\usepackage{bm}
\usepackage{latexsym, amscd ,psfrag}
\usepackage{graphicx}
\usepackage{caption}
\usepackage{float}

\usepackage[mathscr]{euscript}
\usepackage{enumitem}
\usepackage{cleveref}

\makeatletter
\newsavebox{\@brx}
\newcommand{\llangle}[1][]{\savebox{\@brx}{\(\m@th{#1\langle}\)}%
  \mathopen{\copy\@brx\kern-0.5\wd\@brx\usebox{\@brx}}}
\newcommand{\rrangle}[1][]{\savebox{\@brx}{\(\m@th{#1\rangle}\)}%
  \mathclose{\copy\@brx\kern-0.5\wd\@brx\usebox{\@brx}}}
\makeatother

\newcommand{\ep}{\epsilon}

\newcommand{\Ga}{\Gamma}

\newcommand{\bC}{\mathbb{C}}

\newcommand{\bL}{\mathbb{L}}
\newcommand{\bP}{\mathbb{P}}

\newcommand{\bZ}{\mathbb{Z}}

\newcommand{\cB}{\mathcal{B}}

\newcommand{\cD}{\mathcal{D}}

\newcommand{\cL}{\mathcal{L}}
\newcommand{\cI}{\mathcal{I}}
\newcommand{\cM}{\mathcal{M}}
\newcommand{\cO}{\mathcal{O}}

\newcommand{\cT}{\mathcal{T}}
\newcommand{\cU}{\mathcal{U}}

\newcommand{\cX}{\mathcal{X}}
\newcommand{\cY}{\mathcal{Y}}
\newcommand{\cZ}{\mathcal{Z}}

\newcommand{{\inv} }{\mathrm{inv}}
\newcommand{\ev}{\mathrm{ev}}
\newcommand{\Aut}{\mathrm{Aut}}

\newcommand{\val}{ {\mathrm{val}} }
\newcommand{\vir}{{\mathrm{vir}}}

\newcommand{\Mbar}{\overline{\cM}}

\newtheorem{lma}{Lemma}[section]

\newtheorem{prop}[lma]{Proposition}

\newtheorem{thm}[lma]{Theorem}

\newtheorem{remark}[lma]{Remark}

\newtheorem{dummy}{dummy}[section]
\theoremstyle{definition}
\newtheorem{definition}[dummy]{Definition}

\newcommand{\Mb}{\Mbar^\bullet_{\chi, \gamma}(\bP^1_{a}, \mu)}
\newcommand{\MP}{\Mbar^\bullet_\chi(\bP^1,\nu,\mu)}
\newcommand{\MQ}{\Mbar^\bullet_{\chi^1}(\bP^1,\nu,\mu)}
\newcommand{\Md}{\Mbar^\bullet_{\chi^0, \gamma-\nu}(\cB\bZ_{a})}

\newcommand{\xm}[1]{{#1}^\bullet_{\chi,\nu,\mu}}
\newcommand{\xn}[1]{{#1}^\bullet_{\chi,\mu,\gamma}}
\newcommand{\xo}[1]{{#1}^\bullet_{\mu}}

\newcommand{\tpi  }{\tilde{\pi}  }

\newcommand{\amm}{|\Aut(\nu)||\Aut(\mu)|}
\newcommand{\am}{|\Aut(\mu)|}

\DeclareMathOperator{\Gcd}{gcd}

\DeclareMathOperator{\Br}{Br}

\DeclareMathOperator{\rk}{rk}

\begin{document}

\title{On one-leg orbifold topological vertex in refined Gromov-Witten theory}

\author{Jinghao Yu}
\address{Jinghao Yu, Department of Mathematical Sciences, Tsinghua University, Haidian District, Beijing 100084, China}
\email{yjh21@mails.tsinghua.edu.cn}

\author{Zhengyu Zong}
\address{Zhengyu Zong, Department of Mathematical Sciences,
	Tsinghua University, Haidian District, Beijing 100084, China}
\email{zyzong@mail.tsinghua.edu.cn}

\maketitle

\begin{abstract}
	We define the one-leg orbifold topological vertex in refined Gromov-Witten theory \cite{BS24}. There are two cases where the leg is effective or gerby. The main result of this paper is the computation of the effective case. In the smooth case, this result matches the one-leg refined vertex in \cite{IKV09}. As an application, we compute the refined Gromov-Witten invariants of the local football.
\end{abstract}

\tableofcontents

\section{Introduction}

\subsection{Historical background}\label{sec:History}
\subsubsection{Topological vertex and its orbifold generalization}
In \cite{AKMV05}, M. Aganagic, A. Klemm, M. Mari\~{n}o and C. Vafa proposed a theory on computing Gromov-Witten invariants in all genera of any smooth toric Calabi-Yau threefold. The building block for their algorithm is the
\textit{topological vertex}, a generating function of open Gromov-Witten invariants which they predict via large N duality. Mathematical study of the topological vertex can be found in \cite{ORV06},\cite{OP10},\cite{Liu-Liu-Zhou1,Liu-Liu-Zhou2,LLLZ09}

The mathematical theory of the topological vertex can also be viewed as the Gromov-Witten/Donaldson-Thomas correspondence \cite{MNOP06,MNOP06b}, which conjectures that the GW theory and the DT theory of a smooth 3-fold are equivalent after a change of variables. For smooth toric Calabi-Yau 3-folds, the GW theory is obtained by gluing the GW vertex, a generating function of cubic Hodge integrals, and the DT theory is obtained by gluing
the DT vertex, a generating function of 3d partitions. The GW/DT
correspondence for smooth toric Calabi-Yau 3-folds can be reduced to
the correspondence between the GW vertex and the DT vertex. In \cite{MOOP11}, the GW/DT
correspondence is proved for smooth toric 3-folds.

A vertex formalism for the orbifold DT theory (resp. orbifold GW
theory) of toric Calabi-Yau 3-orbifolds is established in \cite{BCY12} (resp. \cite{Ross14}).
For toric Calabi-Yau 3-orbifolds, the orbifold GW theory is obtained by
gluing the GW orbifold vertex, a generating function of cubic abelian
Hurwitz-Hodge integrals, and the orbifold DT theory is obtained by gluing the DT orbifold vertex, a generating function of colored 3d partitions. Mathematical study of the correspondence between the GW orbifold vertex and the DT orbifold vertex can be found in \cite{Zong15,RZ13,RZ15}.

\subsubsection{Refined topological vertex}
In \cite{IKV09}, A. Iqbal, C. Kozçaz and C. Vafa introduced the refined topological vertex, which can be viewed as the building block to study the partition functions of refined topological string on toric Calabi-Yau 3-folds. In \cite{NO16}, Nekrasov-Okounkov study the K-theoretic Donaldson-Thomas/Pandharipande-Thomas (PT) partition functions of 3-folds through an equivariant index of a compactification
of the moduli of immersed curves on an associated local Calabi-Yau fivefold. In the case of toric Calabi-Yau 3-folds, the K-theoretic PT partition function is determined by the K-theoretic vertex. The refined topological vertex can be obtained from the K-theoretic vertex by taking certain limit determined by the preferred slope. The mathematical study of the K-theoretic vertex and the refined topological vertex was further developed in \cite{Arb21,KOO21}.

In \cite{BS24}, Brini-Schuler introduced refined Gromov-Witten theory of a class of Calabi-Yau 3-folds via the equivariant Gromov-Witten theory of the associated Calabi-Yau 5-folds. They introduced the refined GW/PT correspondence, relating the refined GW partition function to the K-theoretic PT partition function in \cite{NO16}. In the case of toric Calabi-Yau 3-folds, this conjecture reduces to the correspondence between refined GW vertex and the K-theoretic vertex.

\subsection{One-leg orbifold refined topological vertex}
Motivated by the works in \ref{sec:History}, we study the one-leg orbifold topological vertex in refined Gromov-Witten theory in this paper. We will define the one-leg vertex in both effective case and gerby case, via refined relative Gromov-Witten theory. We will solve the effective case completely and leave the gerby case in future work.
\subsubsection{The effective case}\label{sec:eff}
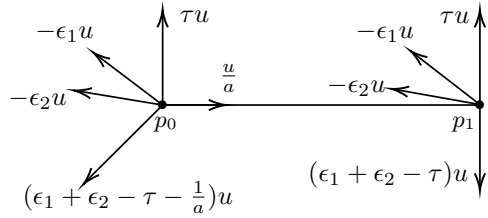
\begin{figure}[H]
\center
\tikzset{every picture/.style={line width=0.65pt}} %set default line width to 0.75pt        

\begin{tikzpicture}[x=0.60pt,y=0.60pt,yscale=-1,xscale=1]
%uncomment if require: \path (0,424); %set diagram left start at 0, and has height of 424

%Straight Lines [id:da7334901322790686] 
\draw    (200,200) -- (400,200) ;
%Straight Lines [id:da7860048751318542] 
\draw    (400,142) -- (400,252) ;
\draw [shift={(400,254)}, rotate = 270] [color={rgb, 255:red, 0; green, 0; blue, 0 }  ][line width=0.75]    (10.93,-3.29) .. controls (6.95,-1.4) and (3.31,-0.3) .. (0,0) .. controls (3.31,0.3) and (6.95,1.4) .. (10.93,3.29)   ;
\draw [shift={(400,140)}, rotate = 90] [color={rgb, 255:red, 0; green, 0; blue, 0 }  ][line width=0.75]    (10.93,-3.29) .. controls (6.95,-1.4) and (3.31,-0.3) .. (0,0) .. controls (3.31,0.3) and (6.95,1.4) .. (10.93,3.29)   ;
%Straight Lines [id:da14673932587529215] 
\draw    (200,200) -- (200,142) ;
\draw [shift={(200,140)}, rotate = 90] [color={rgb, 255:red, 0; green, 0; blue, 0 }  ][line width=0.75]    (10.93,-3.29) .. controls (6.95,-1.4) and (3.31,-0.3) .. (0,0) .. controls (3.31,0.3) and (6.95,1.4) .. (10.93,3.29)   ;
%Straight Lines [id:da2914220884854658] 
\draw    (200,200) -- (151.28,248.46) ;
\draw [shift={(150,250)}, rotate = 312.81] [color={rgb, 255:red, 0; green, 0; blue, 0 }  ][line width=0.75]    (10.93,-3.29) .. controls (6.95,-1.4) and (3.31,-0.3) .. (0,0) .. controls (3.31,0.3) and (6.95,1.4) .. (10.93,3.29)   ;
%Straight Lines [id:da6546981152430535] 
\draw    (200,200) -- (238,200) ;
\draw [shift={(240,200)}, rotate = 180] [color={rgb, 255:red, 0; green, 0; blue, 0 }  ][line width=0.75]    (10.93,-3.29) .. controls (6.95,-1.4) and (3.31,-0.3) .. (0,0) .. controls (3.31,0.3) and (6.95,1.4) .. (10.93,3.29)   ;
%Straight Lines [id:da7153349339200789] 
\draw    (200,200) -- (157.74,167.99) ;
\draw [shift={(156.14,166.79)}, rotate = 37.14] [color={rgb, 255:red, 0; green, 0; blue, 0 }  ][line width=0.75]    (10.93,-3.29) .. controls (6.95,-1.4) and (3.31,-0.3) .. (0,0) .. controls (3.31,0.3) and (6.95,1.4) .. (10.93,3.29)   ;
%Straight Lines [id:da7926484723398272] 
\draw    (200,200) -- (146.12,191.11) ;
\draw [shift={(144.14,190.79)}, rotate = 9.37] [color={rgb, 255:red, 0; green, 0; blue, 0 }  ][line width=0.75]    (10.93,-3.29) .. controls (6.95,-1.4) and (3.31,-0.3) .. (0,0) .. controls (3.31,0.3) and (6.95,1.4) .. (10.93,3.29)   ;
%Straight Lines [id:da9151267347557622] 
\draw    (400,199.7) -- (357.74,166.99) ;
\draw [shift={(356.14,165.79)}, rotate = 37.14] [color={rgb, 255:red, 0; green, 0; blue, 0 }  ][line width=0.75]    (10.93,-3.29) .. controls (6.95,-1.4) and (3.31,-0.3) .. (0,0) .. controls (3.31,0.3) and (6.95,1.4) .. (10.93,3.29)   ;
%Straight Lines [id:da60735066210283] 
\draw    (400,199.7) -- (346.12,190.11) ;
\draw [shift={(344.14,189.79)}, rotate = 9.37] [color={rgb, 255:red, 0; green, 0; blue, 0 }  ][line width=0.75]    (10.93,-3.29) .. controls (6.95,-1.4) and (3.31,-0.3) .. (0,0) .. controls (3.31,0.3) and (6.95,1.4) .. (10.93,3.29)   ;

% Text Node
\draw (209,138.4) node [anchor=north west][inner sep=0.75pt]    {$\tau u$};
% Text Node
\draw (373,139.4) node [anchor=north west][inner sep=0.75pt]    {$\tau u$};
% Text Node
\draw (110,246.4) node [anchor=north west][inner sep=0.75pt]    {$(\epsilon_1+\epsilon_2-\tau-\frac{1}{a})u$};
% Text Node
\draw (290,233.4) node [anchor=north west][inner sep=0.75pt]    {$(\ep_1+\ep_2-\tau)u$};
% Text Node
\draw (234,172.4) node [anchor=north west][inner sep=0.75pt]    {$\frac{u}{a}$};
% Text Node
\draw (118,146.4) node [anchor=north west][inner sep=0.75pt]    {$-\epsilon _{1}u$};
% Text Node
\draw (101,188.4) node [anchor=north west][inner sep=0.75pt]    {$-\epsilon _{2}u$};
% Text Node
\draw (320,146.4) node [anchor=north west][inner sep=0.75pt]    {$-\epsilon _{1}u$};
% Text Node
\draw (302,180.4) node [anchor=north west][inner sep=0.75pt]    {$-\epsilon _{2}u$};
% Text Node
\draw (194.5,195) node [anchor=north west][inner sep=0.75pt]    {\footnotesize $\bullet$};
% Text Node
\draw (394.5,195) node [anchor=north west][inner sep=0.75pt]    {\footnotesize $\bullet$};
% Text Node
\draw (193,205.4) node [anchor=north west][inner sep=0.75pt]    {{\small $p_0$}};
% Text Node
\draw (381,205.4) node [anchor=north west][inner sep=0.75pt]    {{\small $p_1$}};

\end{tikzpicture}
\caption{The weight diagram of $Z$}\label{fig:weights-i}
\end{figure}

Let $a$ be a positive integer. Let $\bP^1_a$ be the projective line $P^1$
with root construction of order $a$ at 0. Let $p_0\in\bP^1_a$ be the unique point with isotropy group $\bZ_a$. The total space $Z:=\textrm{Tot}(\cO_{\bP^1_a}(-p_0)\oplus \cO_{\bP^1_a}\oplus \cO_{\bP^1_a}\oplus \cO_{\bP^1_a}\to \bP^1_a)$ is a toric 5-orbifold and we let $\widetilde{T}\cong(\bC^*)^5\subset Z$ be the 5-dimensional embedded torus in $Z$. Let $p_1\in \bP^1_a$ be the other $\widetilde{T}$-fixed point and let $D_\infty$ be the fiber of $\cO_{\bP^1_a}(-p_0)\oplus \cO_{\bP^1_a}\oplus \cO_{\bP^1_a}\oplus \cO_{\bP^1_a}$ at $p_1$. Then $Z':=Z\setminus D_\infty$ is an affine toric Calabi-Yau 5-orbifold. 
Let $T\cong (\bC^*)^4\subset \widetilde{T}$ be the 4-dimensional Calabi-Yau sub-torus which acts trivially on the canonical bundle of $Z'$. We will also consider a 1-dimensional sub-torus $\bC^*\subset T$ and let $\phi:H^*(BT) \to H^*(B\bC^*)=\bC[u]$ be the induced homomorphism. Let the $\bC^*$ weights at the tangent space $T_{p_0}\bP^1_a$ and at $\cO_{\bP^1_a}(-p_0)|_{p_0}, \ \cO_{\bP^1_a}|_{p_0}, \ \cO_{\bP^1_a}|_{p_0}, \ \cO_{\bP^1_a}|_{p_0}$ be $\frac{u}{a},(\epsilon_1+\epsilon_2-\tau-\frac{1}{a})u,\tau u,-\epsilon_1 u,-\epsilon_2 u$ respectively (See Figure \ref{fig:weights-i}).

In order to define the refined topological vertex, we introduce the {\bf framing condition}:
$$
 \tau=\epsilon_1.
$$
Let $\gamma = (\gamma_1,\dots, \gamma_n)$ be the vectors of integers
$$  
  1\leq \gamma_i \leq a-1.
$$
Let $\mu$ be a partition of $d>0$. In Section \ref{sec:refined GW}, we will define $\overline{\cM}_{g,\gamma}(\bP^1_a,\mu)$ which is the moduli space of relative maps to $(\bP^1_a,\infty)$ with ramification profile $\mu$ over $p_1$. A point in $\overline{\cM}_{g,\gamma}(\bP^1_a,\mu)$ is of the form
$$
[f: (C,x_1,\dots,x_n,y_1,\dots,y_{l(\mu)})\rightarrow (\bP^1_a,p_1^{(m)})]
$$
such that 
\[
f^{-1}(p_1^{(m)}) = \sum_{i=1}^{l(\mu)}\mu_iy_i
\]
as Cartier divisors and the monodromy around $x_i$ is given by $\gamma_i\in \bZ_a$. We will also consider the disconnected version $\Mb$, where the domain curve $C$ is allowed to be disconnected with
$2(h^0(\cO_C)-h^1(\cO_C))=\chi$. Let
$$\pi:\cU\to \Mb$$
be the universal domain curve and let $\cT$ be the universal target. Then there is an evaluation map
$$F:\cU\to \cT$$
and a contraction map
$$\tpi: \cT\to \bP^1_a.$$
Let $\cD\subset \cU$ be the divisor corresponding to the $l(\mu)$ marked points $y_1,\cdots,y_{l(\mu)}$. Define
\begin{eqnarray*}
	V_D&=&V_{\ep_1}=V_{\ep_2}=R^1\pi_*(\tilde{F}^*\cO_{\bP^1_a}(-\mathcal{D}) ),\\
	V_{D_d}&=&R^1\pi_* \tilde{F}^*\cO_{\bP^1_a}(-p_{0}),\\
\end{eqnarray*}
where $\tilde{F}=\tilde\pi\circ F:\mathcal{U}\to \bP^1_a,p_1=\infty\in\bP^1_a$, and $p_{0}$ is the stack point of $\bP^1_a$. The $\bC^*$ action on $Z$ induces $\bC^*$ actions on $\Mb$ and on $	V_D,V_{\ep_1},V_{\ep_2},V_{D_d}$.

Consider the following equivariant relative GW invariant:
\begin{equation*}
		\xn{\tilde{K}} :=\frac{1}{\am|\Aut(\gamma)|}\int_{[\Mb^{\bC^*}]^{\vir}}\frac{e_{\bC^*}(V_D\oplus V_{D_d}\oplus V_{\ep_1}\oplus V_{\ep_2})}{e_{\bC^*}(N^\vir)},
\end{equation*}
where $\overline{\cM}^\bullet_{\chi,\gamma}(\bP^1_a,\mu)^{\bC^*}$ is the fixed locus of the $\bC^*$-action on $\overline{\cM}^\bullet_{\chi,\gamma}(\bP^1_a,\mu)$, $N^\vir$ is the virtual normal bundle of 
$\overline{\cM}^\bullet_{\chi,\gamma}(\bP^1_a,\mu)^{\bC^*}$ in $\overline{\cM}^\bullet_{\chi,\gamma}(\bP^1_a,\mu)$.
We define the \textit{one-leg effective refined vertex} as the following generating function
$$
\tilde{K}^\bullet_\mu(x) := \sqrt{-1}^{-d}(\ep_1u)^{-l(\mu)}\sum_{\chi,\gamma}x_\gamma\xn{\tilde{K}},
$$
where $x_\gamma:=x_{\gamma_1}\cdots x_{\gamma_n}$.

\subsubsection{The gerby case}
Let $m$ be a positive integer. Let $\cY:=\bP^1\times B\bZ_m$ be the trivial $B\bZ_m$ gerbe over $\bP^1$. Let $\cL$ be the $m$-th root of the line bundle $\cO_{\cY}$. In other words, $\cL^{\otimes m}=\cO_{\cY}$ and the group $\bZ_m$ acts on $\cL$ fiberwisely via the fundamental representation. We have a natural $\bC^*$ action on $\cY$ and let $p_0,p_1\in\cY$ be the $\bC^*$ fixed points. 
As in Section \ref{sec:eff}, we consider four line bundles $\cL^{-1}(-1):=\cL^{-1}\otimes \cO_{\cY}(-1),\cL,\cO_{\cY},\cO_{\cY}$ over $\cY$ and let the $\bC^*$ weights at the tangent space $T_{p_0}\cY$ and at $\cL^{-1}\otimes \cO_{\cY}(-1)|_{p_0}$, $\cL|_{p_0},\ \cO_{\cY}|_{p_0},\ \cO_{\cY}|_{p_0}$ be $w_1u,(-w_1-w_2+\epsilon_1+\epsilon_2)u,w_2u,-\epsilon_1 u,-\epsilon_2 u$ respectively (See Figure \ref{fig:weights-ii}), where $H^*(B\bC^*)=\bC[u]$.

\begin{figure}[H]
\center
\tikzset{every picture/.style={line width=0.65pt}} %set default line width to 0.75pt        

\begin{tikzpicture}[x=0.60pt,y=0.60pt,yscale=-1,xscale=1]
%uncomment if require: \path (0,424); %set diagram left start at 0, and has height of 424

%Straight Lines [id:da7334901322790686] 
\draw    (200,200) -- (400,200) ;
%Straight Lines [id:da7860048751318542] 
\draw    (400,142) -- (400,252) ;
\draw [shift={(400,254)}, rotate = 270] [color={rgb, 255:red, 0; green, 0; blue, 0 }  ][line width=0.75]    (10.93,-3.29) .. controls (6.95,-1.4) and (3.31,-0.3) .. (0,0) .. controls (3.31,0.3) and (6.95,1.4) .. (10.93,3.29)   ;
\draw [shift={(400,140)}, rotate = 90] [color={rgb, 255:red, 0; green, 0; blue, 0 }  ][line width=0.75]    (10.93,-3.29) .. controls (6.95,-1.4) and (3.31,-0.3) .. (0,0) .. controls (3.31,0.3) and (6.95,1.4) .. (10.93,3.29)   ;
%Straight Lines [id:da14673932587529215] 
\draw    (200,200) -- (200,142) ;
\draw [shift={(200,140)}, rotate = 90] [color={rgb, 255:red, 0; green, 0; blue, 0 }  ][line width=0.75]    (10.93,-3.29) .. controls (6.95,-1.4) and (3.31,-0.3) .. (0,0) .. controls (3.31,0.3) and (6.95,1.4) .. (10.93,3.29)   ;
%Straight Lines [id:da2914220884854658] 
\draw    (200,200) -- (151.28,248.46) ;
\draw [shift={(150,250)}, rotate = 312.81] [color={rgb, 255:red, 0; green, 0; blue, 0 }  ][line width=0.75]    (10.93,-3.29) .. controls (6.95,-1.4) and (3.31,-0.3) .. (0,0) .. controls (3.31,0.3) and (6.95,1.4) .. (10.93,3.29)   ;
%Straight Lines [id:da6546981152430535] 
\draw    (200,200) -- (238,200) ;
\draw [shift={(240,200)}, rotate = 180] [color={rgb, 255:red, 0; green, 0; blue, 0 }  ][line width=0.75]    (10.93,-3.29) .. controls (6.95,-1.4) and (3.31,-0.3) .. (0,0) .. controls (3.31,0.3) and (6.95,1.4) .. (10.93,3.29)   ;
%Straight Lines [id:da7153349339200789] 
\draw    (200,200) -- (157.74,167.99) ;
\draw [shift={(156.14,166.79)}, rotate = 37.14] [color={rgb, 255:red, 0; green, 0; blue, 0 }  ][line width=0.75]    (10.93,-3.29) .. controls (6.95,-1.4) and (3.31,-0.3) .. (0,0) .. controls (3.31,0.3) and (6.95,1.4) .. (10.93,3.29)   ;
%Straight Lines [id:da7926484723398272] 
\draw    (200,200) -- (146.12,191.11) ;
\draw [shift={(144.14,190.79)}, rotate = 9.37] [color={rgb, 255:red, 0; green, 0; blue, 0 }  ][line width=0.75]    (10.93,-3.29) .. controls (6.95,-1.4) and (3.31,-0.3) .. (0,0) .. controls (3.31,0.3) and (6.95,1.4) .. (10.93,3.29)   ;
%Straight Lines [id:da9151267347557622] 
\draw    (400,199.7) -- (357.74,166.99) ;
\draw [shift={(356.14,165.79)}, rotate = 37.14] [color={rgb, 255:red, 0; green, 0; blue, 0 }  ][line width=0.75]    (10.93,-3.29) .. controls (6.95,-1.4) and (3.31,-0.3) .. (0,0) .. controls (3.31,0.3) and (6.95,1.4) .. (10.93,3.29)   ;
%Straight Lines [id:da60735066210283] 
\draw    (400,199.7) -- (346.12,190.11) ;
\draw [shift={(344.14,189.79)}, rotate = 9.37] [color={rgb, 255:red, 0; green, 0; blue, 0 }  ][line width=0.75]    (10.93,-3.29) .. controls (6.95,-1.4) and (3.31,-0.3) .. (0,0) .. controls (3.31,0.3) and (6.95,1.4) .. (10.93,3.29)   ;

% Text Node
\draw (209,138.4) node [anchor=north west][inner sep=0.75pt]    {$w_2 u$};
% Text Node
\draw (361,139.4) node [anchor=north west][inner sep=0.75pt]    {$w_2 u$};
% Text Node
\draw (96,246.4) node [anchor=north west][inner sep=0.75pt]    {$(-w_1-w_2+\ep_1+\ep_2)u$};
% Text Node
\draw (268,233.4) node [anchor=north west][inner sep=0.75pt]    {$(-w_2+\ep_1+\ep_2)u$};
% Text Node
\draw (220,180.4) node [anchor=north west][inner sep=0.75pt]    {$w_1u$};
% Text Node
\draw (118,148.4) node [anchor=north west][inner sep=0.75pt]    {$-\epsilon _{1}u$};
% Text Node
\draw (101,188.4) node [anchor=north west][inner sep=0.75pt]    {$-\epsilon _{2}u$};
% Text Node
\draw (318,148.4) node [anchor=north west][inner sep=0.75pt]    {$-\epsilon _{1}u$};
% Text Node
\draw (302,180.4) node [anchor=north west][inner sep=0.75pt]    {$-\epsilon _{2}u$};
% Text Node
\draw (194.5,195) node [anchor=north west][inner sep=0.75pt]    {\footnotesize $\bullet$};
% Text Node
\draw (394.5,195) node [anchor=north west][inner sep=0.75pt]    {\footnotesize $\bullet$};
% Text Node
\draw (193,205.4) node [anchor=north west][inner sep=0.75pt]    {{\small $p_0$}};
% Text Node
\draw (381,205.4) node [anchor=north west][inner sep=0.75pt]    {{\small $p_1$}};

\end{tikzpicture}
\caption{The weight diagram of $\cY$}\label{fig:weights-ii}
\end{figure}

We impose the framing condition 
$$
w_1=\epsilon_1.
$$

Let $\gamma = (\gamma_1,\dots, \gamma_n)$ be the vectors of integers
$$
1\leq \gamma_i \leq m-1.
$$ 
Let $\overline{\mu}=\{(\mu_1,k_1),\cdots,(\mu_l,k_l)\}$ be a $\bZ_m$ weighted partition of $d>0$, meaning that $\{\mu_1,\cdots,\mu_l\}$ is a partition of $d$ and $k_i\in\bZ_m,i=1,\cdots,l$. Let $l(\overline{\mu})=l$ denote the length of the $\bZ_m$ weighted partition $\overline{\mu}$. As in Section \ref{sec:eff}, we can similarly consider the moduli space $\overline{\cM}_{g,\gamma}(\cY,\overline{\mu})$ of relative stable maps to $(\cY,p_1)$ with ramification profile $\overline{\mu}$ and the disconnected version $\overline{\cM}^\bullet_{\chi,\gamma}(\cY,\overline{\mu})$. Let
$$\pi:\cU\to\overline{\cM}^\bullet_{\chi,\gamma}(\cY,\overline{\mu})$$
be the universal domain curve and let $\cT$ be the universal target. Then there is an evaluation map
$$F:\cU\to \cT$$
and a contraction map
$$\tpi: \cT\to \cY.$$
Let $\cD\subset \cU$ be the divisor corresponding to the $l(\overline{\mu})$ marked points $y_1,\cdots,y_{l(\overline{\mu})}$. Define
\begin{eqnarray*}
	V_D&=&R^1\pi_*(\tilde{F}^*\cL(-\mathcal{D}) )\\
	V_{\ep_1}&=&V_{\ep_2}=R^1\pi_*(\tilde{F}^*\cO_{\cY}(-\mathcal{D}) ),\\
	V_{D_d}&=&R^1\pi_* \tilde{F}^*\cL^{-1}(-1),\\
\end{eqnarray*}
where $\tilde{F}=\tilde\pi\circ F:\mathcal{U}\to \cY$.
The $\bC^*$ action on $$\cL^{-1}(-1),\cL,\cO_{\cY},\cO_{\cY}$$ induces $\bC^*$ actions on $\overline{\cM}^\bullet_{\chi,\gamma}(\cY,\overline{\mu})$ and on $V_D,V_{\ep_1},V_{\ep_2},V_{D_d}$.

Consider the following equivariant relative GW invariant:
\begin{equation*}
	\tilde{K}^{\bullet,\textrm{gerby}}_{\chi,\overline{\mu},\gamma} :=\frac{1}{|\Aut(\overline{\mu})||\Aut(\gamma)|}\int_{[\overline{\cM}^\bullet_{\chi,\gamma}(\cY,\overline{\mu})^{\bC^*}]^{\vir}}\frac{e_{\bC^*}(V_D\oplus V_{D_d}\oplus V_{\ep_1}\oplus V_{\ep_2})}{e_{\bC^*}(N^\vir)},
\end{equation*}
where $\overline{\cM}^\bullet_{\chi,\gamma}(\cY,\overline{\mu})^{\bC^*}$ is the fixed locus of the $\bC^*$-action on $\overline{\cM}^\bullet_{\chi,\gamma}(\cY,\overline{\mu})$, $N^\vir$ is the virtual normal bundle of 
$\overline{\cM}^\bullet_{\chi,\gamma}(\cY,\overline{\mu})^{\bC^*}$ in $\overline{\cM}^\bullet_{\chi,\gamma}(\cY,\overline{\mu})$.
We define the \textit{one-leg gerby refined vertex} as the following generating function
$$
\tilde{K}^{\bullet,\textrm{gerby}}_{\overline{\mu}}(x) := (\ep_1u)^{-l(\overline{\mu})}\sum_{\chi,\gamma}x_\gamma\tilde{K}^{\bullet,\textrm{gerby}}_{\chi,\overline{\mu},\gamma},
$$
where $x_\gamma:=x_{\gamma_1}\cdots x_{\gamma_n}$.

\subsection{Statement of the main result}
For a partition $\mu$, $\chi_{\mu}$ denotes the character of the irreducible representation of $S_d$
indexed by $\mu$, where $d=|\mu| = \sum_{i=1}^{l(\mu)} \mu_i$. Let $\mu'$ be the conjugate partition of $\mu$. Let $s_{\mu}$ be the Schur function corresponding to $\mu$.

The following theorem is the main result of this paper.

\begin{thm}[= Theorem \ref{thm:refine-vertex}]\label{thm:main}
	Let $t=e^{-\ep_1 u}, q_l=\xi_a^{-1}e^{-\sum_{i=1}^{a-1}\frac{w_a^{-2il}}{a}(w_a^i-w_a^{-i})x_i},l=1,
	\cdots,a-1$, where $\xi_a=e^{\frac{2\pi i}{a}},w_a=e^{\frac{\pi i}{a}}$. Then in the framing condition $\tau = \ep_1$, we have
	\begin{equation*}
		\begin{aligned}
			\xo{\tilde{K}}(x) = \sqrt{-1}^{|\mu|}(t^{\frac{1}{2}}q_1^{-\frac{1}{a}}\cdots q_{a-1}^{-\frac{a-1}{a}})^{|\mu|}
			\sum_{|\nu|=|\mu|}s_{\nu'}(\tilde{t}_\bullet)\frac{\chi_\nu(\mu)}{z_\mu},
		\end{aligned}
	\end{equation*} 
	where $\tilde{t}_{\bullet} = (\tilde{T},\tilde{T}q_{a-1},\dots, \tilde{T}q_1\cdots q_{a-1})$ and $\tilde{T} = (1,t,t^2,t^3,\dots)$. In the smooth case (when $a=1$ and $x = \emptyset$), let $t^{-\rho} = (t^{\frac{1}{2}},t^{\frac{3}{2}},\dots)$, we have
	\[
	\xo{\tilde{K}} = \sqrt{-1}^{|\mu|}\sum_{|\nu|=|\mu|}s_{\nu'}(t^{-\rho})\frac{\chi_{\nu}(\mu)}{z_\mu}. 
	\]
\end{thm}

In the smooth case (when $a=1$), the expression of the one-leg vertex in Theorem \ref{thm:main} is consistent with the result in \cite{IKV09}.

\subsection{Overview of the paper}
In Section \ref{sec:classical-MV}, we review the generalized {Mari\~{n}o-Vafa formula studied in \cite{Zong15}. In Section \ref{sec:refined GW}, we define the one-leg effective refined vertex. In Section \ref{sec:loc}, we use virtual localization to compute the one-leg effective refined vertex. In Section \ref{sec:main result}, we impose the framing condition and obtain the formula for the one-leg effective refined vertex, which is the main result of this paper. In Section \ref{sec:app}, we apply our result to compute the refined Gromov-Witten invariants of the local football.

\subsection*{Acknowledgements}
The authors would like to thank Bohan Fang, Chiu-Chu Melissa Liu, and Song Yu for useful discussions. The second author is partially supported by the Natural Science Foundation of Beijing, China (grant No. 1252008) and NSFC (grant No. 12571067).

\section{Mari\~{n}o-Vafa formula in orbifold Gromov-Witten theory}\label{sec:classical-MV}
In this section, we review the generalized {Mari\~{n}o-Vafa formula studied in \cite{Zong15}. This formula can be viewed as the expression of the effective one-leg orbifold topological vertex in unrefined Gromov-Witten theory.

Fix an integer $a\geq 1$ and let $\cB\bZ_a$ be the classifying stack of principal $\bZ_a$ bundles. Let $\Mbar_{g, \gamma}(\cB\bZ_a)$ be the moduli space of stable maps from genus $g$, $n$-pointed twisted curves to $\cB\bZ_a$, where $\gamma=(\gamma_{1}, \cdots, \gamma_{n})$ is a vector of elements in $\bZ_a$ indicating the monodromy data around the marked points. Let $U$ be the irreducible representation of $\bZ_a$ defined as
$$\phi^U:\bZ_a\to \bC^*,\phi^U(1)=e^{\frac{2\pi\sqrt{-1}}{a}}.$$
Then there is a corresponding Hodge bundle
$$E^U\to\Mbar_{g, \gamma}(\cB\bZ_a)$$
and the corresponding Hodge classes on $\Mbar_{g, \gamma}(\cB\bZ_a)$ are defined by Chern classes of $E^U$,
$$\lambda^U_i=c_i(E^U).$$
More generally, given any irreducible representation $R$ of $\bZ_a$, we have a corresponding Hodge bundle $E^R$ and Hodge classes $\lambda^R_i$. Let $\Mbar_{g,n}$ be the moduli space of genus $g$, $n$-pointed stable curves and let $\psi_i$ be the $i^{th}$ descendant class on $\Mbar_{g,n}$, $1\leq i\leq n$. Define
$$\epsilon :\Mbar_{g, \gamma}(\cB\bZ_a)\to \Mbar_{g,n}$$
to be the canonical morphism which forgets the orbifold structures. We define the descendant classes $\bar{\psi}_i$ on $\Mbar_{g, \gamma}(\cB\bZ_a)$ as
$$\bar{\psi}_i=\epsilon^*(\psi_i).$$
Given a formal variable $u$, we define
$$
\Lambda^{\vee,R}_g(u):=u^{\rk E^R}-\lambda_1^R u^{\rk E^R-1} +\cdots+(-1)^{\rk E^R}\lambda_{\rk E^R}^R,
$$
where $\rk E^R$ is the rank of $E^R$ determined by the orbifold Riemann-Roch formula.

Let $d$ be a positive integer and let $\mu=(\mu_{1}\geq\cdots\geq\mu_{l(\mu)}>0)$ be a partition of $d>0$ meaning that $|\mu|:=\sum_{i=1}^{l(\mu)}\mu_i=d$. Now we require $\gamma=(\gamma_{1}, \cdots, \gamma_{n})$ to be a vector of \emph{nontrivial} elements in $\bZ_a$. We view $\mu$ as a vector of elements in $\bZ_a$ by considering $\mu_i$ modulo $a\bZ$ for each $\mu_i$. Then for $\tau\in \bZ$, we define $G_{g,\mu,\gamma}(\tau)_{a}$ as
\begin{eqnarray*}
G_{g,\mu,\gamma}(\tau)_{a}&:=&\frac{\sqrt{-1}^{l(\mu)-|\mu|+2\sum_{i=1}^{l(\mu)}[\frac{\mu_i}{a}]}
\tau ^{l(\mu)-1}a^{l(\mu)-\sum_{i=1}^{l(\mu)}\delta_{0,\langle \frac{\mu_i}{a}\rangle}}}{\am|\Aut(\gamma)|}
\prod_{i=1}^{l(\mu)}\frac{\prod_{l=1}^{[\frac{\mu_i}{a}]}(\mu_i\tau+l)}{[\frac{\mu_i}{a}]!}
\\
&&\cdot\int_{\Mbar_{g, \gamma-\mu}(\cB\bZ_a)}\frac{\left(-\frac{1}{a}(\tau+\frac{1}{a})\right)^{-\delta}
\Lambda_{g}^{\vee,U}(\frac{1}{a})\Lambda_{g}^{\vee,U^\vee}(-\tau-\frac{1}{a})
\Lambda_{g}^{\vee,1}(\tau)}{\prod_{i=1}^{l(\mu)}(1-\mu_i\bar{\psi}_i)},
\end{eqnarray*}
where $\gamma-\mu$ denotes the vector $(\gamma_{1}, \cdots, \gamma_{n},-\mu_1,\cdots,-\mu_{l(\mu)})$, $\bar{\psi}_i$ corresponds to $\mu_i$, $U^\vee$ and 1 denote the dual of $U$ and the trivial representation respectively, $[x]$ denotes the integer part of $x$, $\langle x\rangle=x-[x]$,
$$\delta_{0,x}=\left\{\begin{array}{ll}1, &x=0,\\
0, &x\neq 0,\end{array} \right.$$
and
$$\delta=\left\{\begin{array}{ll}1, &\textrm{if all monodromies around loops on the domain curve are trivial,}\\
0, &\textrm{otherwise}.\end{array} \right.$$
Introduce formal variables $p=(p_1,p_2,\ldots,p_n,\ldots),x=(x_1,\ldots,x_{a-1})$ and define
$$
p_\mu=p_{\mu_1}\cdots p_{\mu_{l(\mu)} },x_\gamma=x_{\gamma_1}\cdots x_{\gamma_n},
$$
where $\mu$ is a partition and $\gamma=(\gamma_{1}, \cdots, \gamma_{n})$ is a vector of nontrivial elements in $\bZ_a$. Define
generating functions
\begin{eqnarray*}
G_{\mu,\gamma}(\hbar;\tau)_{a}&=&\sum_{g=0}^{\infty}\hbar^{2g-2+l(\mu)}G_{g,\mu,\gamma}(\tau)_{a}\\
G(\hbar;\tau;p;x)_{a}&=&\sum_{\mu\neq \emptyset,\gamma}G_{\mu,\gamma}(\hbar;\tau)_{a}p_\mu x_\gamma=\sum_{\mu\neq \emptyset}G_{\mu}(\hbar;\tau;x)_{a}p_\mu=\sum_{\gamma}G_{\gamma}(\hbar;\tau;p)_{a}x_\gamma\\
G^\bullet(\hbar;\tau;p;x)_{a}&=&\exp(G(\hbar;\tau;p;x)_{a})=
\sum_{\mu,\gamma}G^\bullet_{\mu,\gamma}(\hbar;\tau)_{a}p_\mu x_\gamma=
1+\sum_{\mu\neq\emptyset,\gamma}G^\bullet_{\mu,\gamma}(\hbar;\tau)_{a}p_\mu x_\gamma\\
&=&1+\sum_{\mu\neq\emptyset}G^\bullet_{\mu}(\hbar;\tau;x)_{a}p_\mu\\
G^\bullet_{\mu,\gamma}(\hbar;\tau)_{a}&=&\sum_{\chi\in 2\bZ,\chi\leq 2l(\mu)}\hbar^{-\chi+l(\mu)}G^\bullet_{\chi,\mu,\gamma}(\tau)_{a}.
\end{eqnarray*}

The generating function $G(\hbar;\tau;p;x)_{a}$ corresponds to the framed 1-leg orbifold Gromov-Witten vertex where the leg is effective \cite{Zong15}. 

For a partition $\mu$, $\chi_{\mu}$ denotes the character of the irreducible representation of $S_d$
indexed by $\mu$, where $d=|\mu| = \sum_{i=1}^{l(\mu)} \mu_i$. Let
$$\kappa_{\mu} = |\mu| + \sum_i (\mu_i^2 - 2i\mu_i).$$

Introduce new variables $q_0,q_1,\cdots,q_{a-1}$ and consider the change of variables
$$
q_0=-e^{\sqrt{-1}\hbar},q_l=\xi_a^{-1}e^{-\sum_{i=1}^{a-1}\frac{w_a^{-2il}}{a}(w_a^i-w_a^{-i})x_i},l=1,
\cdots,a-1,
$$
where $\xi_a=e^{\frac{2\pi i}{a}},w_a=e^{\frac{\pi i}{a}}$.

Then under the above change of variables, let
\begin{eqnarray*}
R^\bullet_{\mu}(\hbar;\tau;x)_{a}&=&\sum_{|\nu|=|\mu|}(q_0^{\frac{1}{2}}q_1^{-\frac{1}{a}}\cdots q_{a-1}^{-\frac{a-1}{a}})^{|\nu|}\sum_{|\xi|=|\nu|}s_{\xi'}(-q_{\bullet})\chi_{\xi}(\nu)
\Phi^\bullet_{\nu,\mu}(\sqrt{-1}\hbar\tau)
\end{eqnarray*}
where
\begin{eqnarray*}
\Phi^\bullet_{\nu,\mu}(\hbar)=\sum_{\eta}\frac{\chi_{\eta}(\nu)}{z_{\nu}} \frac{\chi_{\eta}(\mu)}{z_{\mu}}e^{\frac{\kappa_{\eta}\hbar}{2}},
\end{eqnarray*}
$z_\mu = |\Aut(\mu)|z_{\mu_1}\dots z_{\mu_{l(\mu)}}$.
$s_{\xi'}$ is the Schur polynomial corresponding to the conjugate of $\xi$, $-q_{\bullet}=(-Q,-Qq_{a-1},\cdots, -Qq_1\cdots q_{a-1})$, and $-Q=(1,-q_0,(-q_0)^2,(-q_0)^3,\cdots)$.

The generating function $\Phi^\bullet_{\nu,\mu}(\hbar)$ satisfies the following properties (see \cite{Liu-Liu-Zhou2}),

\begin{equation}\label{eqn:Comp}
  \Phi^{\bullet}_{\nu, \mu}(\hbar_1+\hbar_2)
=  \sum_{\sigma} \Phi^{\bullet}_{\nu, \sigma}(\hbar_1)
\cdot z_{\sigma} \cdot \Phi^{\bullet}_{\sigma, \mu}(\hbar_2),
\end{equation}
\begin{equation}\label{eqn:Init}
 \Phi^{\bullet}_{\nu, \mu}(0) = \frac{1}{z_{\nu}}\delta_{\nu, \mu}.
\end{equation}
By (\ref{eqn:Comp}), (\ref{eqn:Init}) we have
\begin{eqnarray*}
R^\bullet_{\mu}(\hbar;0;x)_{a}&=&(q_0^{\frac{1}{2}}q_1^{-\frac{1}{a}}\cdots q_{a-1}^{-\frac{a-1}{a}})^{|\mu|}\sum_{|\nu|=|\mu|}s_{\nu'}(-q_{\bullet})\frac{\chi_{\nu}(\mu)}{z_\mu},
\end{eqnarray*}
\begin{equation}\label{eqn:R}
R^\bullet_{\mu}(\hbar;\tau;x)_{a}=\sum_{|\nu|=|\mu|}R^\bullet_{\nu}(\hbar;0;x)_{a}z_{\nu}
\Phi^\bullet_{\nu,\mu}(\sqrt{-1}\hbar\tau).
\end{equation}

We also define generating functions

\begin{eqnarray*}
R^\bullet(\hbar;\tau;p;x)_{a}&=&\sum_{\mu,\gamma}R^\bullet_{\mu}(\hbar;\tau;x)_{a}p_\mu,\\
R(\hbar;\tau;p;x)_{a}&=&\log R^\bullet(\hbar;\tau;p;x)_{a}=\sum_{\mu\neq\emptyset,\gamma}R_{\mu,\gamma}(\hbar;\tau)_{a}p_\mu x_\gamma\\
&=&\sum_{\mu\neq \emptyset}R_{\mu}(\hbar;\tau;x)_{a}p_\mu=\sum_{\gamma}R_{\gamma}(\hbar;\tau;p)_{a}x_\gamma.
\end{eqnarray*}

The generating series $G^\bullet_\mu(\hbar;\tau;x)_a$ is determined by its initial value in the following way.
\begin{thm}[\cite{Zong15}]\label{key} It holds that
\begin{eqnarray*}
G^\bullet_{\mu}(\hbar;\tau;x)_{a}=\sum_{|\nu|=|\mu|}G^\bullet_{\nu}(\hbar;0;x)_{a}z_{\nu}
\Phi^\bullet_{\nu,\mu}(\sqrt{-1}\tau\hbar).
\end{eqnarray*}
Equivalently, $$G^\bullet_{\mu}(\hbar;0;x)_{a}=\sum_{|\nu|=|\mu|}G^\bullet_{\nu}(\hbar;\tau;x)_{a}z_{\nu}
\Phi^\bullet_{\nu,\mu}(-\sqrt{-1}\tau\hbar).$$
\end{thm}
Theorem \ref{key} and \eqref{eqn:R} imply the Mari\~{n}o-Vafa formula for $\bZ_a$:
\begin{thm}[Generalized Mari\~{n}o-Vafa formula \cite{Zong15}]\label{MVza}
\begin{eqnarray*}
G(\hbar;\tau;p;x)_{a}=R(\hbar;\tau;p;x)_{a}.
\end{eqnarray*}
\end{thm}

When $a=1$, Theorem \ref{MVza} is the original Mari\~{n}o-Vafa formula proved in \cite{Liu-Liu-Zhou1} and \cite{Oko-Pan}.

\section{One-leg orbifold topological vertex in refined Gromov-Witten theory}\label{sec:refined GW}
In this section, we define the effective one-leg orbifold topological vertex in refined Gromov-Witten theory. The readers are referred to \cite{BS24} for a study of general refined Gromov-Witten theory.
\subsection{Moduli spaces and $\bC^*$-actions}
Let $\bP^1_a$ be the projective line $\bP^1$ with root construction of order $a$ at 0.
For an integer $m>0$, let 
\[
  \bP^1_a[m] =\bP^1_a\cup \bP^1_{(1)}\cup \dots \cup \bP^1_{(m)}
\]
be the union of $\bP^1_a$ and a chain of $m$ copies $\bP^1$, where $\bP^1_a$ is glued to $\bP^1_{(1)}$ at $p_1^{(0)}$
and $\bP^1_{(l)}$ is glued to $\bP^1_{(l+1)}$ at $p_1^{(l)}$ for $1\leq l\leq m-1$. We call the irreducible component $\bP^1_a$
the root component and the other irreducible components the bubble components. Denote by $\pi[m]: \bP^1_a[m]\rightarrow \bP^1_a$ the map 
which is identity on the root component and contracts all the bubble components to $p^{(0)}_1$.
Let 
\[
  \bP^1(m) = \bP^1_{(1)}\cup \dots \cup \bP^1_{(m)}
\]
denote the union of bubble components of $\bP^1_a[m]$. For convenience, we set $\bP^1_a[0]=\bP^1_{(0)}=\bP^1_a$.

Let $\gamma = (\gamma_1,\dots, \gamma_n)$ be the vectors of integers
\[
  1\leq \gamma_i \leq a-1.
\]
Let $\mu$ be a partition of $d>0$. Define $\overline{\cM}_{g,\gamma}(\bP^1_a,\mu)$ to be the moduli space of relative stable maps to $(\bP^1_a,\infty)$.
Then a point in $\overline{\cM}_{g,\gamma}(\bP^1_a,\mu)$ is of the form
$$
  [f: (C,x_1,\dots,x_n,y_1,\dots,y_{l(\mu)})\rightarrow (\bP^1_a,p_1^{(m)})]
$$
where $C$ is a genus $g$ prestable curve and
\[
    f^{-1}(p_1^{(m)}) = \sum_{i=1}^{l(\mu)}\mu_iy_i
\]
as Cartier divisors. We require the parity condition
\[
    d=\sum_{i=1}^{n}\gamma_i \pmod a
\]
so that the moduli space $\overline{\cM}_{g,\gamma}(\bP^1_a,\mu)$ is non-empty.
We will also consider the disconnected version $\Mb$, where the domain curve $C$ is allowed to be disconnected with
$2(h^0(\cO_C)-h^1(\cO_C))=\chi$.
\par
Similarly, if we specify ramification types $\nu,\mu$ over $0,\infty\in\bP^1$, we can define the corresponding moduli space $\overline{\cM}_{g,0}(\bP^1,\nu,\mu)$ and 
$\overline{\cM}^\bullet_\chi(\bP^1,\nu,\mu)$ of relative stable maps.

Consider the $\bC^*$-action
\[
  t \cdot [z^0:z^1]=[tz^0:z^1]
\]
on $\bP^1$, where $[z^0:z^1]$ is the homogeneous coordinates. This action lifts canonically on $\bP^1_a$. These induce actions on $\bP^1[m]$ and on $\bP^1_a[m]$ with 
trivial actions on the bubble components. These in turn induce actions on $\overline{\cM}_{g,\gamma}(\bP^1_a,\mu), \overline{\cM}^\bullet_{\chi,\gamma}(\bP^1_a,\mu), \overline{\cM}_{g,0}(\bP^1,\nu,\mu)$, and 
$\overline{\cM}^\bullet_\chi(\bP^1,\nu,\mu)$. Define the quotient space $\overline{\cM}^\bullet_\chi(\bP^1,\nu,\mu)//\bC^*$ to be
$$
  \overline{\cM}^\bullet_\chi(\bP^1,\nu,\mu)//\bC^* = (\overline{\cM}^\bullet_\chi(\bP^1,\nu,\mu)\backslash \overline{\cM}^\bullet_\chi(\bP^1,\nu,\mu)^{\bC^*})/\bC^*.
$$

\subsection{The branch morphism and double Hurwitz numbers}
Let $\bP^1[m_0,m_1]$ be a chain of $\bP^1$ obtained by attaching $\bP^1(m_0)$ and $\bP^1(m_1)$ to $\bP^1$ at $0$ and $\infty$ respectively.
Let $q_{m_0}^0$ and $q_{m_1}^1$ be the distinguished point on $\bP^1[m_0,m_1]$. 
Similar to the case of $\overline{\cM}_{g,\gamma}(\bP^1_a,\mu)$, a map $[f]\in\overline{\cM}^\bullet_\chi(\bP^1,\nu,\mu)$ has target at $\bP^1[m_0,m_1]$,
and the ramification profiles of $[f]$ at $q_{m_0}^0$ and $q_{m_1}^1$ are $\nu$ and $\mu$ respectively.
Let $\pi[m_0,m_1]:\bP^1[m_0,m_1]\rightarrow\bP^1$ be the contraction to the root
component. Let $r=-\chi+l(\nu)+l(\mu)$ be the virtual dimension of $\overline{\cM}^\bullet_\chi(\bP^1,\nu,\mu)$. Then there is a branch morphism
$$ \textrm{Br}:\Mbar^\bullet_{\chi}(\bP^1,\nu, \mu)\to \textrm{Sym}^r\bP^1$$
sending $[f:C\to \bP^1[m_0,m_1]]$ to
$$\textrm{div} (\tilde{f})-(d-l(\nu))0-(d-l(\mu))\infty,$$
where div$(\tilde{f})$ is the branch divisor of $\tilde{f}=\pi [m_{0},m_1]\circ f:C\to \bP^1$.

The disconnected double Hurwitz number is defined as
$$
\xm{H} =\frac{1}{\amm}\int_{[\MP]^{\vir} }\Br^*( H^{r}).
$$
By virtual localization, we have the following proposition \cite{Liu-Liu-Zhou2}. 
\begin{prop}[\cite{Liu-Liu-Zhou2}]
$$\xm{H}=\frac{r!}{\amm}\int_{[\MP//\bC^*]^{\vir} }(\psi^0)^{r-1},$$
where $\psi^0$ is the target $\psi$ class, the first Chern class of the line bundle $\bL_0$ over $\MP$ whose fiber at
$$[f:C\to \bP^1[m_{0},m_1]]$$
is the cotangent line $T^*_{q^0_{m_0}}\bP^1[m_{0},m_1]$.
\end{prop}

The double Hurwitz number $\xm{H}$ satisfies the following Burnside formula \cite{Dij}:
\begin{prop}[Burnside formula \cite{Dij}]
$$\Phi^\bullet_{\nu,\mu}(\hbar)=\sum_{\chi}\hbar^{-\chi+l(\mu)+l(\nu)}\frac{\xm{H}}{(-\chi+l(\mu)+l(\nu))!}$$
where $\Phi^\bullet_{\nu,\mu}(\hbar)$ is defined combinatorially in Section \ref{sec:classical-MV}.
\end{prop}
As a corollary, it holds that
\begin{equation}\label{eqn:minus-phi}
  \Phi^\bullet_{\nu,\mu}(-\hbar) = (-1)^{l(\mu)+l(\nu)}\Phi^\bullet_{\nu,\mu}(\hbar).
\end{equation}

\subsection{The obstruction bundle}
Let
$$\pi:\cU\to \Mb$$
be the universal domain curve and let $\cT$ be the universal target. Then there is an evaluation map
$$F:\cU\to \cT$$
and a contraction map
$$\tpi: \cT\to \bP^1_a.$$
Let $\cD\subset \cU$ be the divisor corresponding to the $l(\mu)$ marked points $y_1,\cdots,y_{l(\mu)}$. Define
\begin{eqnarray*}
V_D&=&R^1\pi_*(\cO_{\cU}(-\mathcal{D}) )\\
V_{D_d}&=&R^1\pi_* \tilde{F}^*\cO_{\bP^1_a}(-p_{0}),
\end{eqnarray*}
where $\tilde{F}=\tilde\pi\circ F:\mathcal{U}\to \bP^1_a,p_1=\infty\in\bP^1_a$, and $p_{0}$ is the stack point of $\bP^1_a$. The fibers of $V_D$ and $V_{D_d}$ at
$$
\left[ f:(C,x_1,\ldots,x_n,y_1,\cdots,y_{l(\mu)})\to \bP^1_a[m]\ \right]\in \Mb
$$
are $H^1(C, \cO_C(-D))$ and $H^1(C, \tilde{f}^*\cO_{\bP^1_a}(-p_0))$,
respectively, where $D=y_1+\ldots+y_{l(\mu)}$ and $\tilde{f}=\pi[m]\circ f$. It is easy to see that the rank of the obstruction bundle
$$V=V_D\oplus V_{D_{d}}$$
is equal to the virtual dimension of $\Mb$, which is $-\chi+n+l(\mu)+\frac{d}{a}-\sum_{i=1}^{n}\frac{\gamma_i}{a}$.

We lift the $\bC^*$-action to the obstruction bundle $V$ as follows. It suffices to lift the $\bC^*$-action on $\bP^1_a$ to the line bundles $\cO_{\bP^1_a}(-p_{0})$ and $\cO_{\bP^1_a}$. Let the weights of the $\bC^*$-action on $\cO_{\bP^1_a}(-p_{0})$ at $p_0$ and $p_1$ be $-w-\frac{1}{a}$ and $-w$, respectively, and let the weights of the $\bC^*$-action on $\cO_{\bP^1_a}$ at $p_0$ and $p_1$ be $\tau$ and $\tau$, respectively, where $\tau,w\in\bZ$. In other words, if we write the obstruction bundle in the form of equivariant divisors, we have
\begin{eqnarray*}
V_D&=&R^1\pi_*(\tilde{F}^*\cO_{\bP^1_a}(\tau(ap_0-p_1))(-\mathcal{D}) )\\
V_{D_d}&=&R^1\pi_* \tilde{F}^*\cO_{\bP^1_a}(-p_{0}+w(p_1-ap_0)).
\end{eqnarray*}

\subsection{Refined relative Gromov-Witten invariants}
Let 
$$V_{\ep_i} = R^1\pi_*(\tilde{F}^*\cO_{\bP^1_a}(-\ep_i(ap_0-p_1))(-\cD)), \quad i = 1,2.$$
The rank of the bundle $V_{\ep_i}$ over $\Mb$ is $l(\mu)-\chi/2$. The construction of $V_{\ep_i}$ is the same as the bundle $V_D$, which is obtained by lifting the $\bC^*$-action on $\bP^1_a$ to the bundle $\cO_{\bP^1_a}$ by weights $-\ep_i$. 

Let $Z$ be the toric 5-fold defined by the total space 
\[
  Z = \text{Tot}(V_D \oplus V_{D_d}\oplus V_{\ep_1}\oplus V_{\ep_2}) \rightarrow \bP^1_a.
\]
The weights of the $\bC^*$ action at $0$ are $\frac{u}{a},\tau u,(-w-\frac{1}{a})u,-\ep_1 u,-\ep_2 u$ respectively (See Figure \ref{fig:weights}).
We impose the following \textit{Calabi-Yau condition} of the $\bC^*$-action:
\[
    -\ep_1 -\ep_2 +\tau -w = 0.
\]

Consider the following equivariant relative Gromov-Witten invariant:
\begin{equation*}
  \begin{aligned}
    \xn{\tilde{K}} &:=\frac{1}{\am|\Aut(\gamma)|}\int_{[\Mb^{\bC^*}]^{\vir}}\frac{e_{\bC^*}(V_D\oplus V_{D_d}\oplus V_{\ep_1}\oplus V_{\ep_2})}{e_{\bC^*}(N^\vir)},   
  \end{aligned}
\end{equation*}
where $\overline{\cM}^\bullet_{\chi,\gamma}(\bP^1_a,\mu)^{\bC^*}$ is the fixed locus of the $\bC^*$-action on $\overline{\cM}^\bullet_{\chi,\gamma}(\bP^1_a,\mu)$, $N^\vir$ is the virtual normal bundle of 
$\overline{\cM}^\bullet_{\chi,\gamma}(\bP^1_a,\mu)^{\bC^*}$ in $\overline{\cM}^\bullet_{\chi,\gamma}(\bP^1_a,\mu)$.

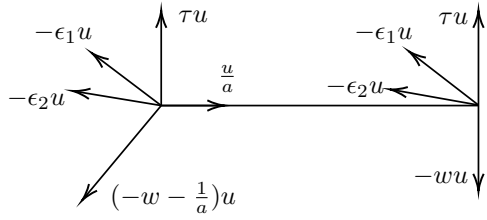
\begin{figure}[H]
\tikzset{every picture/.style={line width=0.65pt}} %set default line width to 0.75pt        
\begin{tikzpicture}[x=0.60pt,y=0.60pt,yscale=-1,xscale=1]
%uncomment if require: \path (0,424); %set diagram left start at 0, and has height of 424

%Straight Lines [id:da7334901322790686] 
\draw    (200,200) -- (400,200) ;
%Straight Lines [id:da7860048751318542] 
\draw    (400,142) -- (400,252) ;
\draw [shift={(400,254)}, rotate = 270] [color={rgb, 255:red, 0; green, 0; blue, 0 }  ][line width=0.75]    (10.93,-3.29) .. controls (6.95,-1.4) and (3.31,-0.3) .. (0,0) .. controls (3.31,0.3) and (6.95,1.4) .. (10.93,3.29)   ;
\draw [shift={(400,140)}, rotate = 90] [color={rgb, 255:red, 0; green, 0; blue, 0 }  ][line width=0.75]    (10.93,-3.29) .. controls (6.95,-1.4) and (3.31,-0.3) .. (0,0) .. controls (3.31,0.3) and (6.95,1.4) .. (10.93,3.29)   ;
%Straight Lines [id:da14673932587529215] 
\draw    (200,200) -- (200,142) ;
\draw [shift={(200,140)}, rotate = 90] [color={rgb, 255:red, 0; green, 0; blue, 0 }  ][line width=0.75]    (10.93,-3.29) .. controls (6.95,-1.4) and (3.31,-0.3) .. (0,0) .. controls (3.31,0.3) and (6.95,1.4) .. (10.93,3.29)   ;
%Straight Lines [id:da2914220884854658] 
\draw    (200,200) -- (151.28,258.46) ;
\draw [shift={(150,260)}, rotate = 309.81] [color={rgb, 255:red, 0; green, 0; blue, 0 }  ][line width=0.75]    (10.93,-3.29) .. controls (6.95,-1.4) and (3.31,-0.3) .. (0,0) .. controls (3.31,0.3) and (6.95,1.4) .. (10.93,3.29)   ;
%Straight Lines [id:da6546981152430535] 
\draw    (200,200) -- (238,200) ;
\draw [shift={(240,200)}, rotate = 180] [color={rgb, 255:red, 0; green, 0; blue, 0 }  ][line width=0.75]    (10.93,-3.29) .. controls (6.95,-1.4) and (3.31,-0.3) .. (0,0) .. controls (3.31,0.3) and (6.95,1.4) .. (10.93,3.29)   ;
%Straight Lines [id:da7153349339200789] 
\draw    (200,200) -- (157.74,167.99) ;
\draw [shift={(156.14,166.79)}, rotate = 37.14] [color={rgb, 255:red, 0; green, 0; blue, 0 }  ][line width=0.75]    (10.93,-3.29) .. controls (6.95,-1.4) and (3.31,-0.3) .. (0,0) .. controls (3.31,0.3) and (6.95,1.4) .. (10.93,3.29)   ;
%Straight Lines [id:da7926484723398272] 
\draw    (200,200) -- (146.12,191.11) ;
\draw [shift={(144.14,190.79)}, rotate = 9.37] [color={rgb, 255:red, 0; green, 0; blue, 0 }  ][line width=0.75]    (10.93,-3.29) .. controls (6.95,-1.4) and (3.31,-0.3) .. (0,0) .. controls (3.31,0.3) and (6.95,1.4) .. (10.93,3.29)   ;
%Straight Lines [id:da9151267347557622] 
\draw    (400,199.7) -- (357.74,166.99) ;
\draw [shift={(356.14,165.79)}, rotate = 37.14] [color={rgb, 255:red, 0; green, 0; blue, 0 }  ][line width=0.75]    (10.93,-3.29) .. controls (6.95,-1.4) and (3.31,-0.3) .. (0,0) .. controls (3.31,0.3) and (6.95,1.4) .. (10.93,3.29)   ;
%Straight Lines [id:da60735066210283] 
\draw    (400,199.7) -- (346.12,190.11) ;
\draw [shift={(344.14,189.79)}, rotate = 9.37] [color={rgb, 255:red, 0; green, 0; blue, 0 }  ][line width=0.75]    (10.93,-3.29) .. controls (6.95,-1.4) and (3.31,-0.3) .. (0,0) .. controls (3.31,0.3) and (6.95,1.4) .. (10.93,3.29)   ;

% Text Node
\draw (209,138.4) node [anchor=north west][inner sep=0.75pt]    {$\tau u$};
% Text Node
\draw (373,139.4) node [anchor=north west][inner sep=0.75pt]    {$\tau u$};
% Text Node
\draw (166,246.4) node [anchor=north west][inner sep=0.75pt]    {$(-w-\frac{1}{a})u$};
% Text Node
\draw (357,237.4) node [anchor=north west][inner sep=0.75pt]    {$-wu$};
% Text Node
\draw (234,172.4) node [anchor=north west][inner sep=0.75pt]    {$\frac{u}{a}$};
% Text Node
\draw (118,146.4) node [anchor=north west][inner sep=0.75pt]    {$-\epsilon _{1}u$};
% Text Node
\draw (101,188.4) node [anchor=north west][inner sep=0.75pt]    {$-\epsilon _{2}u$};
% Text Node
\draw (320,146.4) node [anchor=north west][inner sep=0.75pt]    {$-\epsilon _{1}u$};
% Text Node
\draw (302,180.4) node [anchor=north west][inner sep=0.75pt]    {$-\epsilon _{2}u$};

\end{tikzpicture}
\caption{The weight diagram of $Z$}\label{fig:weights}
\end{figure}

We will calculate $\xn{\tilde{K}}$ in the next section by virtual localization. Define the generating function $\tilde{K}^\bullet_\mu(x)$ to be
$$
  \tilde{K}^\bullet_\mu(x) = \sqrt{-1}^{-d}(\ep_1u)^{-l(\mu)}\sum_{\chi,\gamma}x_\gamma\xn{\tilde{K}}.
$$
We will impose the framing condition $\tau = \ep_1$ in Section \ref{sec:main result} and we call $\tilde{K}^\bullet_\mu(x)|_{\tau = \ep_1}$ the one-leg effective refined vertex.

\section{Virtual localization}\label{sec:loc}
In this section, we calculate the relative version $\xn{\tilde{K}}$ by virtual localization.

\subsection{Fixed points}
The connected components of the $\bC^*$ fixed points set of $\Mb$ are parameterized by labeled graphs. We first introduce some graph notations.

Let
$$\left[ f:(C,x_1,\ldots,x_n,y_1,\cdots,y_{l(\mu)})\to \bP^1_a[m]\right]\in\Mb $$
be a fixed point of the $\bC^*$-action. The restriction of the map
$$ \tilde{f}=\pi[m]\circ f: C\to \bP^1_a$$
to an irreducible component of $C$ is either a constant map to one of the $\bC^*$
fixed points $p_0, p_1$ or a cover of $\bP^1_a$ which is fully
ramified over $p_0$ and $p_1$. We associate a labeled graph $\Gamma$ to
the $\bC^*$ fixed point
$$\left[ f:(C,x_1,\ldots,x_n,y_1,\cdots,y_{l(\mu)})\to \bP^1_a[m] \right]$$ as follows:
\begin{enumerate}
\item We assign a vertex $v$ to each connected
component $C_v$ of $\tilde{f}^{-1}(\{p_0,p_1\})$, a label
$i(v)=i$ if $\tilde{f}(C_v)=p_i$, where $i=0,1$, and a label $g(v)$
which is the arithmetic genus of $C_v$ (we define $g(v)=0$ if $C_v$ is a
point). For $i(v)=0$, we define $n(v)$ to be the number of marked points on $C_v$. Denote by $V(\Gamma)^{(i)}$ the set of vertices with $i(v)=i$,
where $i=0,1$. Then the set $V(\Gamma)$ of vertices of the graph $\Gamma$
is a disjoint union of $V(\Gamma)^{(0)}$ and $V(\Gamma)^{(1)}$. For $v\in V(\Gamma)^{(0)}$ define
\begin{eqnarray*}
r_0(v)=&2g(v)-2 + \val(v)+n(v), & v\in V(\Gamma)^{(0)},
\end{eqnarray*}
\begin{eqnarray*}
V^I(\Gamma)^{(0)}&=&\{v\in V(\Gamma)^{(0)}:g(v)=0,\val(v)=1,n(v)=0\},\\
V^{I,I}(\Gamma)^{(0)}&=&\{v\in V(\Gamma)^{(0)}:g(v)=0,\val(v)=1,n(v)=1\},\\
V^{II}(\Gamma)^{(0)}&=&\{v\in V(\Gamma)^{(0)}:g(v)=0,\val(v)=2,n(v)=0\},\\
V^S(\Gamma)^{(0)}&=&\{v\in V(\Gamma)^{(0)}:r_0(v)> 0\}.
\end{eqnarray*}

\item We assign an edge $e$ to each rational irreducible component $C_e$ of $C$
such that $\tilde{f}|_{C_e}$ is not a constant map.
Let $d(e)$ be the degree of $\tilde{f}|_{C_e}$.
Then $\tilde{f}|_{C_e}$ is fully ramified over $p_0$ and $p_1$.
Let $E(\Gamma)$ denote the set of edges of $\Gamma$.

\item The set of flags of $\Gamma$ is given by
$$F(\Gamma)=\{(v,e):v\in V(\Gamma), e\in E(\Gamma),
C_v\cap C_e\neq \emptyset \}.$$

\item For each $v\in V(\Gamma)$, define
$$d(v)=\sum_{(v,e)\in F(\Gamma)}d(e),$$
and let $\nu(v)$ be the partition of $d(v)$
determined by $\{d(e): (v,e)\in F(\Gamma)\}$ and let $\nu$ be the partition of $d$ determined by $\{d(e): e\in E(\Gamma)\}$.
When the target is $\bP^1_a[m]$, where $m>0$,
we assign an additional label for each
$v\in V(\Gamma)^{(1)}$:
let $\mu(v)$ be the partition of
$d(v)$ determined by the ramification of
$f|_{C_v}:C_v\to\bP^1_a(m)$ over $p_1^{(m)}$.
\end{enumerate}
Note that for $v\in V(\Gamma)^{(1)}$,
$\nu(v)$ coincides with the partition of $d(v)$
determined by the ramification of $f|_{C_v}:C_v\to \bP^1_a(m)$
over $p_1^{(0)}$.

Let $\cM_{\nu_i}$ be the moduli space of $\bC^*$-fixed degree $\nu_i$ covers of $\bP^1_a$ with stack structure given by $\nu_i\textrm{ (mod $a$)}$. Then the $\bC^*$-fixed locus can be identified with
$$\bigsqcup_{\chi^0+\chi^1-2l(v)=\chi}(\Md\times_{\bar{I}\cB\bZ_{a}^{l(\nu)}}\cM_{\nu_1}\times\cdots\times\cM_{\nu_{l(\nu
)}}\times\MQ//\bC^*)/\Aut(\nu),$$
where $\bar{I}\cB\bZ_{a}$ is the rigidified inertia stack of $\cB\bZ_{a}$. Therefore, we can calculate our integral over
$$\bigsqcup_{\chi^0+\chi^1-2l(v)=\chi}\Md\times\MQ//\bC^*$$
provided we include the following factor:

$$\frac{1}{|\Aut(\nu)|}\prod_{i=1}^{l(\nu)}\frac{1}{\nu_i}\frac{a}{b_i},$$
where $b_i=\frac{a}{\textrm{gcd}(a,\nu_i)}$ is the order of $\nu_i\in \bZ_a$.

\subsection{Virtual normal bundle}
\subsubsection{The target is $\bP^1_a$}
Introduce 
$$\delta_{v}=\left\{\begin{array}{ll}1, &\textrm{if all monodromies around loops on $C_v$ are trivial},\\
0, &\textrm{otherwise}.\end{array} \right.$$
We have the following Feynman rules:
$$\frac{1}{e_{\bC^*}(N_{\Gamma}^{\mathrm{vir} })}
= \prod_{v \in V(\Gamma)^{(0)}} A_v
\prod_{e \in E(\Gamma)} A_e,
$$
where
\begin{eqnarray*}
A_v &= &\left\{
\begin{array}{ll}
\Lambda_{g(v)}^{\vee,U}(\frac{u}{a})(\frac{u}{a})^{-\delta_v}\prod_{(v,e)\in F(\Gamma)}
\frac{(\frac{u}{a})^{\delta_{0,\langle \frac{d(e)}{a}\rangle}}}{(\frac{u}{d(e)} - \bar{\psi}_{(v, e)})\frac{\Gcd(a,d(e))}{a}},
   & v \in V^S(\Gamma)^{(0)}, \\
\frac{u}{d(v)}, & v \in V^I(\Gamma)^{(0)}, \\
1, & v \in V^{I,I}(\Gamma)^{(0)},\\
\frac{1}{\frac{\Gcd(a,d(e_1))}{d(e_1)} + \frac{\Gcd(a,d(e_2))}{d(e_2)}}(\frac{u}{a})^{-1+\delta_{0,\langle \frac{d(e_1)}{a}\rangle}},
   & v \in V^{II}(\Gamma)^{(0)},\\
&(v, e_1), (v, e_2) \in F(\Gamma),
\end{array} \right.\\
 A_e&=& \frac{d(e)^{[\frac{d(e)}{a}]}}{[\frac{d(e)}{a}]!} u^{-[\frac{d(e)}{a}]}.
\end{eqnarray*}

\subsubsection{The target is $\bP^1_a[m]$, $m>0$}
We have the following Feynman rules:
$$\frac{1}{e_{\bC^*}(N_{\Gamma}^{\mathrm{vir} })}
= \frac{1}{-u-\psi^0}\prod_{v \in V(\Gamma)} A_v
\prod_{e \in E(\Gamma)} A_e,
$$
where
\begin{eqnarray*}
A_v &= &\left\{
\begin{array}{ll}
 \prod_{(v,e)\in F(\Gamma)} d(e), &
 v \in V(\Gamma)^{(1)}\\
\Lambda_{g(v)}^{\vee,U}(\frac{u}{a})(\frac{u}{a})^{-\delta_v}\prod_{(v,e)\in F(\Gamma)}
\frac{(\frac{u}{a})^{\delta_{0,\langle \frac{d(e)}{a}\rangle}}}{(\frac{u}{d(e)} - \bar{\psi}_{(v, e)})\frac{\Gcd(a,d(e))}{a}},
   & v \in V^S(\Gamma)^{(0)}, \\
\frac{u}{d(v)}, & v \in V^I(\Gamma)^{(0)}, \\
1, & v \in V^{I,I}(\Gamma)^{(0)},\\
\frac{1}{\frac{\Gcd(a,d(e_1))}{d(e_1)} + \frac{\Gcd(a,d(e_2))}{d(e_2)}}(\frac{u}{a})^{-1+\delta_{0,\langle \frac{d(e_1)}{a}\rangle}},
   & v \in V^{II}(\Gamma)^{(0)},\\
&(v, e_1), (v, e_2) \in F(\Gamma),
\end{array} \right.\\
 A_e&=& \frac{d(e)^{[\frac{d(e)}{a}]}}{[\frac{d(e)}{a}]!} u^{-[\frac{d(e)}{a}]}.
\end{eqnarray*}

\subsection{The bundle $V_{D}$}\label{D}
\subsubsection{The target is $\bP^1_a$}
We have
$$
i_{\Gamma}^* e_{\bC^*}(V_D)=\prod_{v\in V(\Gamma)^{(0)}}A^D_{v},
$$
where
$$
A_{v}^D = \Lambda_{g(v)}^{\vee,1}(\tau u) \cdot
(\tau u)^{\val(v)-1}.
$$

\subsubsection{The target is $\bP^1_a[m]$, $m>0$}
We have the following Feynman rules:
$$i_{\Gamma}^*e_{\bC^*}(V_D) = \prod_{v\in V(\Gamma)} A_v^D,$$
where
$$
 A_v^{D} = \begin{cases}
\Lambda_{g(v)}^{\vee,1}(\tau u) \cdot (\tau u)^{\val(v)-1},
& v \in V(\Gamma)^{(0)}, \\
\Lambda_{g(v)}^{\vee,1}(\tau u) \cdot
(\tau u)^{l(\mu(v))-1},
& v \in V(\Gamma)^{(1)}.
\end{cases} \\
$$

\subsection{The bundle $V_{D_d}$} \label{Dd}
\subsubsection{The target is $\bP^1_a$}
We have the following Feynman rules:
$$i_{\Gamma}^*e_{\bC^*}(V_{D_d}) = \prod_{v\in V(\Gamma)^{(0)}} A_v^{D_d} \cdot
\prod_{e\in E(\Gamma)} A_e^{D_d},$$
where
\begin{eqnarray*}
A_v &= &\left\{
\begin{array}{ll}
\Lambda_{g(v)}^{\vee,U^\vee}((-w-\frac{1}{a})u)((-w-\frac{1}{a})u)^{-\delta_v}\\ \prod_{(v,e)\in F(\Gamma)}
((-w-\frac{1}{a})u)^{\delta_{0,\langle \frac{d(e)}{a}\rangle}},
    &v \in V^S(\Gamma)^{(0)}, \\
1, & v \in V^I(\Gamma)^{(0)}, \\
1, & v \in V^{I,I}(\Gamma)^{(0)},\\
\left((-w-\frac{1}{a})u\right)^{\delta_{0,\langle \frac{d(e_1)}{a}\rangle}},
    &v \in V^{II}(\Gamma)^{(0)},(v, e_1),(v,e_2) \in F(\Gamma),
\end{array} \right.\\
 A_e&=& \frac{\prod_{l=1}^{[\frac{d(e)}{a}]-\delta_{0,\langle \frac{d(e)}{a}\rangle}}(d(e)w+l)}{d(e)^{[\frac{d(e)}{a}]-\delta_{0,\langle \frac{d(e)}{a}\rangle}}} (-u)^{[\frac{d(e)}{a}]-\delta_{0,\langle \frac{d(e)}{a}\rangle}}.
\end{eqnarray*}

\subsubsection{The target is $\bP^1_a[m]$, $m>0$}
We have the following Feynman rules:
$$i_{\Gamma}^*e_{\bC^*}(V_{D_d}) = \prod_{v\in V(\Ga)} A_v^{D_d} \cdot
\prod_{e\in E(\Gamma)} A_e^{D_d},$$
where
\begin{eqnarray*}
A_v &= &\left\{
\begin{array}{ll}
\Lambda_{g(v)}^{\vee,U^\vee}((-w-\frac{1}{a})u)((-w-\frac{1}{a})u)^{-\delta_v}\\ \prod_{(v,e)\in F(\Gamma)}
((-w-\frac{1}{a})u)^{\delta_{0,\langle \frac{d(e)}{a}\rangle}},
    &v \in V^S(\Gamma)^{(0)}, \\
1, & v \in V^I(\Gamma)^{(0)}, \\
1, & v \in V^{I,I}(\Gamma)^{(0)},\\
\left((-w-\frac{1}{a})u\right)^{\delta_{0,\langle \frac{d(e_1)}{a}\rangle}},
    &v \in V^{II}(\Gamma)^{(0)},(v, e_1),(v,e_2) \in F(\Gamma),\\
    \Lambda_{g(v)}^{\vee,1}(-w u) \cdot
(-w u)^{\val(v)-1}, & v \in V(\Gamma)^{(1)},
\end{array} \right.\\
 A_e&=& \frac{\prod_{l=1}^{[\frac{d(e)}{a}]-\delta_{0,\langle \frac{d(e)}{a}\rangle}}(d(e)w+l)}{d(e)^{[\frac{d(e)}{a}]-\delta_{0,\langle \frac{d(e)}{a}\rangle}}} (-u)^{[\frac{d(e)}{a}]-\delta_{0,\langle \frac{d(e)}{a}\rangle}}.
\end{eqnarray*}

\subsection{The bundle $V_{\ep_i}$}
The localization computation of $V_{\ep_i}$ is same as the bundle $V_D$ by substituting $\tau$ by $-\ep_i$. 
\subsubsection{The target is $\bP^1_a$}
$$
i_{\Gamma}^* e_{\bC^*}(V_{\ep_i})=\prod_{v\in V(\Gamma)^{(0)}}A^{\ep_i}_{v},
$$
where
$$
A_{v}^{\ep_i} = \Lambda_{g(v)}^{\vee,1}(-\ep_i u) \cdot
(-\ep_i u)^{\val(v)-1}.
$$

\subsubsection{The target is $\bP^1_a[m]$}
$$i_{\Gamma}^*e_{\bC^*}(V_{\ep_i}) = \prod_{v\in V(\Gamma)} A_v^{\ep_i},$$
where
$$
 A_v^{\ep_i} = \begin{cases}
\Lambda_{g(v)}^{\vee,1}(-\ep_i u) \cdot (-\ep_i u)^{\val(v)-1},
& v \in V(\Gamma)^{(0)}, \\
\Lambda_{g(v)}^{\vee,1}(-\ep_i u) \cdot
(-\ep_i u)^{l(\mu(v))-1},
& v \in V(\Gamma)^{(1)}.
\end{cases} \\
$$

\subsection{Contribution from each graph}\label{sec:cont}
Define
\[
  C_\Ga := \frac{i^*_\Ga (e_{\bC^*}(V_D)e_{\bC^*}(V_{D_d})e_{\bC^*}(V_{\ep_1})e_{\bC^*}(V_{\ep_2}))}{e_{\bC^*}(N_\Ga^\vir)}.
\]
Let $a_\nu:= \nu_1\cdots\nu_{l(\nu)}$. We write $C_\Ga = \tilde{A}^0\tilde{A}^1$ in the following way.
\begin{equation*}
  \begin{aligned}
    \tilde{A}^0 = & \ a_\nu\prod_{i=1}^{l(\nu)}\frac{\prod_{l=1}^{[\frac{\nu_i}{a}]}(\nu_i w+l)}{[\frac{\nu_i}{a}]!(u-\nu_i\bar{\psi}_i)\frac{\gcd(\nu_i,a)}{a}}
    \\
    &\cdot \prod_{v\in V(\Ga)^{(0)}} (-1)^{\sum_{(v,e)\in F(\Ga)}[\frac{d(e)}{a}]-\delta_v}\Lambda^{\vee,U}_{g(v)}(\frac{u}{a})\cdot(\frac{u}{a})^{\sum_{(v,e)\in F(\Ga)}\delta_{0,\langle\frac{d(e)}{a}\rangle}-\delta_v}
    \\
    & \cdot \prod_{v\in V(\Ga)^{(0)}}\Lambda^{\vee,1}_{g(v)}(\tau u)\cdot (\tau u)^{\val(v)-1}
    \\
    & \cdot \prod_{v\in V(\Ga)^{(0)}}\Lambda^{\vee,U^\vee}_{g(v)}((-w-\frac{1}{a})u)\cdot \Big((w+\frac{1}{a})u\Big)^{-\delta_v}
    \\
    & \cdot \prod_{v\in V(\Ga)^{(0)}}\Lambda^{\vee,1}_{g(v)}(-\ep_1 u)\cdot (-\ep_1 u)^{\val(v)-1}
    \cdot \Lambda^{\vee,1}_{g(v)}(-\ep_2 u)\cdot (-\ep_2 u)^{\val(v)-1}.
  \end{aligned}
\end{equation*}
If the target is $\bP^1_a$, then $\tilde{A}^1=1$. If the target is $\bP^1_a[m]$, $m>0$, then
\begin{equation*}
  \begin{aligned}
    \tilde{A}^1 = & \ \frac{a_\nu}{-u-\psi^0}\cdot  \prod_{v\in V(\Ga)^{(1)}}\Lambda_{g(v)}^{\vee,1}(\tau u) \cdot (\tau u)^{l(\mu(v))-1}\cdot\Lambda_{g(v)}^{\vee,1}(-w u) \cdot (-w u)^{\val(v)-1}
    \\
    &\cdot  \prod_{v\in V(\Ga)^{(1)}}\Lambda_{g(v)}^{\vee,1}(-\ep_1 u) \cdot (-\ep_1 u)^{l(\mu(v))-1}\cdot\Lambda_{g(v)}^{\vee,1}(-\ep_2 u) \cdot (-\ep_2 u)^{l(\mu(v))-1}.
  \end{aligned}
\end{equation*}

\section{Vertex formula}\label{sec:main result}

\subsection{Framing condition}
Consider the following restriction of weights for $\bC^*$: 
$$\tau = \ep_1$$
We call the above restriction as the {\bf framing condition}.
In this setting, we can apply the Mumford's relation:
\[
    \Lambda^{\vee,1}_{g(v)}(-\tau u)\Lambda^{\vee,1}_{g(v)}(\tau u) = (-1)^{g(v)}(\tau u)^{2 g(v)}, \ v \in V(\Ga).
\]
In the framing condition, we have
\[
    C_\Ga = A^0A^1, 
\]
where
\begin{eqnarray*}
A^0&=&a_{\nu}\sqrt{-1}^{l(\nu)-|\nu|}(\sqrt{-1}\ep_1 u)^{-\chi^0+l(\nu)+l(\mu)}\prod_{i=1}^{l(\nu)}\frac{\prod_{l=1}^{[\frac{\nu_i}{a}]}(\nu_{i}w+l)}{[\frac{\nu_i}
{a}]!(u-\nu_{i}\bar{\psi_i})\frac{\Gcd(\nu_{i},a)}{a}}\\
&&\cdot\prod_{v \in V(\Gamma)^{(0)}}(-1)^{\sum_{(v,e)\in F(\Gamma)}[\frac{d(e)}{a}]-\delta_v}\Lambda_{g(v)}^{\vee,U}(\frac{u}{a})\Lambda_{g(v)}^{\vee,U^\vee}((-w-\frac{1}{a})u)
\Lambda_{g(v)}^{\vee,1}(w u)\\
&&\cdot\left(\frac{u}{a}\right)^{\sum_{(v,e)\in F(\Gamma)}\delta_{0,\langle \frac{d(e)}{a}\rangle}-\delta_v}\left((w+\frac{1}{a})u\right)^{-\delta_v}(w u)^{\val(v)-1}\\
A^1&=& \left\{\begin{array}{ll}\sqrt{-1}^{|\mu|-l(\mu)}, &\textrm{the target is }\bP^1_a,\\
\sqrt{-1}^{|\nu|-l(\mu)}a_{\nu}\frac{(\sqrt{-1}w u)^{-\chi^{1}+l(\mu)+l(\nu)}(\sqrt{-1}\ep_1 u)^{-\chi^1+l(\mu)+l(\nu)}}{-u-\psi^0},  &\textrm{the target is }\bP^1_a[m],m>0.
\end{array}\right.
\end{eqnarray*}

\subsection{Refined vertex formula in the framing condition}
Under the framing condition
\begin{eqnarray*}
&&\xn{\tilde{K}}\\
&=&\frac{1}{\am|\Aut(\gamma)|}\int_{[\Mb^{\bC^*}]^{\vir}}\frac{e_{\bC^*}(V_D\oplus V_{D_d}\oplus V_{\ep_1}\oplus V_{\ep_2})}{e_{\bC^*}(N^\vir)}\\
&=&\frac{1}{\am|\Aut(\gamma)|}\sum_{\chi^0+\chi^1-2l(\nu)=\chi}\frac{1}{|\Aut(\nu)|}\prod_{i=1}^{l(\nu)}\frac{1}{\nu_i}
\frac{a}{\frac{a}{\textrm{gcd}(a,\nu_i)}}\\
&&\cdot\int_{[\Md\times\MQ//\bC^*]^{\vir}}\frac{i_{\Gamma}^*(e_{\bC^*}(V_D)e_{\bC^*}(V_{D_d})e_{\bC^*}(V_{\ep_1})e_{\bC^*}(V_{\ep_2}))}{e_{\bC^*}(N_{\Gamma}^{\mathrm{vir} })}\\
&=&\sqrt{-1}^{|\nu|-l(\mu)}\left(\sum_{\chi^0+\chi^1-2l(\nu)=\chi}(\sqrt{-1}\ep_1 u)^{-\chi^0+l(\nu)+l(\mu)}G^\bullet_{\chi^0,\nu,\gamma}(w)_{a}
\cdot z_{\nu}
\frac{(w\ep_1u)^{-\chi^1 +l(\nu) +l(\mu)}}
{(-\chi^1 +l(\nu) +l(\mu))!}H^\bullet_{\chi^1,\nu,\mu}\right),
\end{eqnarray*}
where $w=-\ep_1 -\ep_2 +\tau=-\ep_2$.

Recall that
$$\xo{\tilde{K}}(x)=\sqrt{-1}^{-d}(\ep_1u)^{-l(\mu)}\sum_{\chi,\gamma}x_\gamma\xn{\tilde{K}}.$$
We have
\begin{equation*}
    \xo{\tilde{K}}(x)=\sum_{|\nu|=|\mu|}G^\bullet_{\nu}(\sqrt{-1}\ep_1 u;w;x)_{a}z_{\nu}
\Phi^\bullet_{\nu,\mu}(w\ep_1 u).
\end{equation*}
By Theorem \ref{key}, we have 
\[
  \xo{\tilde{K}}(x) = G^\bullet_{\mu}(\sqrt{-1}\ep_1 u; 0;x)_a.
\]
We see that only $\ep_1$ appears in $\xo{\tilde{K}}(x)$. Furthermore, by Theorem \ref{MVza}, we have the following theorem.
\begin{thm}\label{thm:refine-vertex}
  Let $t=e^{-\ep_1 u}, q_l=\xi_a^{-1}e^{-\sum_{i=1}^{a-1}\frac{w_a^{-2il}}{a}(w_a^i-w_a^{-i})x_i},l=1,
\cdots,a-1$, then in the framing condition $\tau = \ep_1$, we have
  \begin{equation*}
    \begin{aligned}
      \xo{\tilde{K}}(x) = \sqrt{-1}^{|\mu|}(t^{\frac{1}{2}}q_1^{-\frac{1}{a}}\cdots q_{a-1}^{-\frac{a-1}{a}})^{|\mu|}
      \sum_{|\nu|=|\mu|}s_{\nu'}(\tilde{t}_\bullet)\frac{\chi_\nu(\mu)}{z_\mu},
    \end{aligned}
  \end{equation*} 
  where $\tilde{t}_{\bullet} = (\tilde{T},\tilde{T}q_{a-1},\dots, \tilde{T}q_1\cdots q_{a-1})$ and $\tilde{T} = (1,t,t^2,t^3,\dots)$. In the smooth case (when $a=1$ and $x = \emptyset$), let $t^{-\rho} = (t^{\frac{1}{2}},t^{\frac{3}{2}},\dots)$, we have
  \[
    \xo{\tilde{K}} = \sqrt{-1}^{|\mu|}\sum_{|\nu|=|\mu|}s_{\nu'}(t^{-\rho})\frac{\chi_{\nu}(\mu)}{z_\mu}. 
\]
\end{thm}

\subsection{Unrefined limit}
In the unrefined limit, 
we set $\ep_1=-\ep_2=\ep$. Then we have
\begin{eqnarray*}
A^0&=&a_{\nu}\sqrt{-1}^{l(\nu)-|\nu|}(\sqrt{-1}\ep u)^{-\chi^0+l(\nu)+l(\mu)}\prod_{i=1}^{l(\nu)}\frac{\prod_{l=1}^{[\frac{\nu_i}{a}]}(\nu_{i}\tau+l)}{[\frac{\nu_i}
{a}]!(u-\nu_{i}\bar{\psi_i})\frac{\Gcd(\nu_{i},a)}{a}}\\
&&\cdot\prod_{v \in V(\Gamma)^{(0)}}(-1)^{\sum_{(v,e)\in F(\Gamma)}[\frac{d(e)}{a}]-\delta_v}\Lambda_{g(v)}^{\vee,U}(\frac{u}{a})\Lambda_{g(v)}^{\vee,U^\vee}((-\tau-\frac{1}{a})u)
\Lambda_{g(v)}^{\vee,1}(\tau u)\\
&&\cdot\left(\frac{u}{a}\right)^{\sum_{(v,e)\in F(\Gamma)}\delta_{0,\langle \frac{d(e)}{a}\rangle}-\delta_v}\left((\tau+\frac{1}{a})u\right)^{-\delta_v}(\tau u)^{\val(v)-1}\\
A^1&=& \left\{\begin{array}{ll}\sqrt{-1}^{|\mu|-l(\mu)}, &\textrm{the target is }\bP^1_a,\\
\sqrt{-1}^{|\nu|-l(\mu)}a_{\nu}\frac{(\sqrt{-1}\tau u)^{-\chi^{1}+l(\mu)+l(\nu)}(\sqrt{-1}\ep u)^{-\chi^1+l(\mu)+l(\nu)}}{-u-\psi^0},  &\textrm{the target is }\bP^1_a[m],m>0.
\end{array}\right.
\end{eqnarray*}
It follows that
$$\xo{\tilde{K}}(x)=\sum_{|\nu|=|\mu|}G^\bullet_{\nu}(\sqrt{-1}\ep u;\tau;x)_{a}z_{\nu}
\Phi^\bullet_{\nu,\mu}(\tau\ep u) = G^\bullet_\mu(\sqrt{-1}\ep u;0;x).$$

\begin{remark}
  $\ep$ plays the same role as $\hbar$ in Theorem \ref{key}. 
\end{remark}

\subsection{Refined topological vertex}
The refined topological vertex \cite{IKV09} is an extension of the ordinary topological vertex that incorporates two refinement parameters, 
$t=e^{\ep_1}$ and $q=e^{-\ep_2}$. It offers a combinatorial framework for computing refined Gromov-Witten/Donaldson-Thomas invariants of toric Calabi-Yau threefolds. 
In its definition, one must choose a preferred direction (usually marked by $||$ in the toric diagram), while the parameters 
$\ep_1,\ep_2$ are assigned to the remaining two legs of the vertex (see Figure \ref{fig:refine-vertex}).

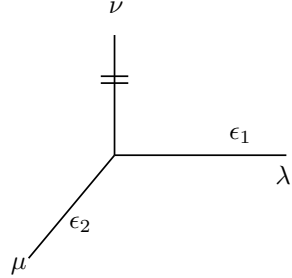
\begin{figure}[H]
	\tikzset{every picture/.style={line width=0.65pt}} %set default line width to 0.75pt        
	\begin{tikzpicture}[x=0.65pt,y=0.65pt,yscale=-1,xscale=1]
		%uncomment if require: \path (0,424); %set diagram left start at 0, and has height of 424
		
		%Straight Lines [id:da315251456999762] 
		\draw    (200,200) -- (300,200) ;
		%Straight Lines [id:da9627600136572271] 
		\draw    (200,200) -- (200,130) ;
		%Straight Lines [id:da4518528549784918] 
		\draw    (200,200) -- (150,260) ;
		%Straight Lines [id:da9498215992746669] 
		\draw    (192,154) -- (208,154) ;
		%Straight Lines [id:da04595534316324368] 
		\draw    (192,158) -- (208,158) ;
		
		% Text Node
		\draw (172,234.4) node [anchor=north west][inner sep=0.75pt]    {$\epsilon_{2}$};
		% Text Node
		\draw (265,182.4) node [anchor=north west][inner sep=0.75pt]    {$\epsilon_{1}$};
		% Text Node
		\draw (292,204.4) node [anchor=north west][inner sep=0.75pt]    {$\lambda $};
		% Text Node
		\draw (138,257.92) node [anchor=north west][inner sep=0.75pt]    {$\mu $};
		% Text Node
		\draw (195,108.92) node [anchor=north west][inner sep=0.75pt]    {$\nu $};

	\end{tikzpicture}
	\caption{refined topological vertex}\label{fig:refine-vertex}
\end{figure}

Given the partitions $\lambda,\mu,\nu$, the refined topological vertex is defined in \cite{IKV09} by
\begin{equation*}
	\begin{aligned}
		C_{\lambda\mu\nu}(t,q) = & \ q^{\frac{\|\mu\|^2+\|\nu\|^2}{2}}t^{-\frac{\|\mu'\|^2}{2}}\tilde{Z}_\nu(t,q)
		\\
		& \times\sum_\eta \Big(\frac{q}{t}\Big)^{\frac{|\eta|+|\lambda|-|\mu|}{2}}s_{\lambda'/\eta}(t^{-\rho}q^{-\nu})s_{\mu/\eta}(q^{-\rho}t^{-\nu'}).
	\end{aligned}
\end{equation*}
where
\[
\tilde{Z}_\nu(t,q) = \prod_{(i,j)\in\nu} (1 - q^{\nu_i - j}t^{\nu_j'-i+1})^{-1}.
\]
In the 1-leg case:
\begin{align*}
	C_{\lambda\emptyset\emptyset}(t,q) &=  \Big(\frac{q}{t}\Big)^{\frac{|\lambda|}{2}}s_{\lambda'}(t^{-\rho}),
	\\
	C_{\emptyset\mu\emptyset}(t,q) &= q^{\frac{\|\mu\|^2}{2}}t^{-\frac{\|\mu'\|^2}{2}}\Big(\frac{q}{t}\Big)^{-\frac{|\mu|}{2}}s_\mu(q^{-\rho}),
	\\
	C_{\emptyset\emptyset\nu}(t,q) &= q^{\frac{\|\nu\|^2}{2}}\tilde{Z}_\nu(t,q).
\end{align*}
The refined smooth GW vertex (i.e. when $a=1$) in Theorem \ref{thm:refine-vertex} is consistent with the first 1-leg refined topological vertex $C_{\lambda\emptyset\emptyset}(t,q)$.

\section{Applications}\label{sec:app}
\subsection{Local football}
Let $\bP^1_{a,b}$ be the football defined by root construction. More concretely, $\bP^1_{a,b}$ is obtained by applying the root construction of order $a,b$ at $0,\infty\in\bP^1$ respectively. Let $p_0$ and $p_\infty$ be the points on $\bP^1_{a,b}$ with isotropy groups $\bZ_a,\bZ_b$ respectively.
Let $\cX_{a,b}$ be the local football defined by
\[
    \cX_{a,b} = \text{Tot}(\cO_{\bP^1_{a,b}}(-p_0)\oplus \cO_{\bP^1_{a,b}}(-p_\infty)).
\]
Let $\cZ_{a,b} = \cX_{a,b}\times \bC^2$, and $T\cong(\bC^*)^4\subset\cZ_{a,b}$ be the Calabi-Yau torus which acts trivially on the canonical bundle of $\cZ_{a,b}$. Let $\bC^*\subset T$ be the one dimensional sub-torus whose weights at $p_0$ and $p_\infty$ are given by the following table:
\begin{equation*}
   \renewcommand{\arraystretch}{1.15}
    \begin{array}{l|cc}
            &   p_0   &  p_\infty \\ [2pt]
            \hline 
            \\ [-11pt]
        T\bP^1_{a,b}    & u/a   & -u/b     \\       
        \cO(-p_0) & (-w-\frac{1}{a})u & -wu\\
        \cO(-p_\infty) & \tau u &  (\tau+\frac{1}{b})u \\
        \cO & -\ep_1 u & -\ep_1 u\\
        \cO & -\ep_2 u & -\ep_2 u
    \end{array}
\end{equation*}
\par
Let $\overline{\cM}_{g,n}(\cZ_{a,b},d)$ be the moduli stack of genus $g$, $n$-pointed stable maps to $\cZ_{a,b}$ of degree $d\in H_2(\cZ_{a,b},\bZ)$. Let $\cI\cZ_{a,b}$ be the inertia stack of $\cZ_{a,b}$,
where
$$\cI\cZ_{a,b} = \bigsqcup_{v\in\bZ\atop {-b+1\leq v\leq a-1}} \cZ_v,$$
and
\begin{equation*}
    \cZ_v = \left\{\begin{aligned}
        &\cB\bZ_a\times \bC^3, && 1 \leq v\leq a-1,
        \\
        &\cZ_{a,b}, &&v=0,
        \\
        & \cB\bZ_b \times \bC^3, && -b+1\leq v \leq -1.
    \end{aligned}\right.
\end{equation*}
Let $\ev_i:\overline{\cM}_{g,n}(\cZ_{a,b},d)\rightarrow \cI\cZ_{a,b}$ be the evaluation map at the $i$-th marked points. Let $\vec{k} =(k_1,\dots,k_{a-1})$, $\vec{l}=(l_1,\dots,l_{b-1})$ be vectors with $k_i,l_j$ as non-negative integers.
Then we associate vectors $\vec{\alpha} = (1^{k_1},2^{k_2},\dots, (a-1)^{k_{a-1}})$ and $\vec{\beta}= (1^{l_1},2^{l_2},\dots,(b-1)^{l_{b-1}})$.
Given any $\vec{\alpha}, \vec{\beta}$, define
\[
    \overline{\cM}_{g,\vec{\alpha},\vec{\beta}}(\cZ_{a,b},d) := \bigcap_{i=1}^{l(\vec{\alpha})}\ev_i^{-1}(\cZ_{\alpha_i})\cap\bigcap_{j=1}^{l(\vec{\beta})}\ev_{l(\vec{\alpha})+j}^{-1}(\cZ_{-\beta_j}).
\]
The $\bC^*$-action on $\cZ_{a,b}$ induces a $\bC^*$-action on $\overline{\cM}_{g,\vec{\alpha},\vec{\beta}}(\cZ_{a,b},d)$. Let $\overline{\cM}_{g,\vec{\alpha},\vec{\beta}}(\cZ_{a,b},d)^{\bC^*}$
denote the $\bC^*$-fixed locus of $\overline{\cM}_{g,\vec{\alpha},\vec{\beta}}(\cZ_{a,b},d)$. The map $\ev_i$ is $\bC^*$-equivariant. 
Define the refined Gromov-Witten invariants of $\cX_{a,b}$ as
\[
    N_{g,d,\vec{\alpha},\vec{\beta}}^{\ep_1,\ep_2} := \int_{[\overline{\cM}_{g,\vec{\alpha},\vec{\beta}}(\cZ_{a,b},d)^{\bC^*}]^{\vir}}\frac{1}{e_{\bC^*}(N^{\vir})},
\]
where $N^\vir$ is the virtual normal bundle of $\overline{\cM}_{g,\vec{\alpha},\vec{\beta}}(\cZ_{a,b},d)^{\bC^*}$ in $\overline{\cM}_{g,\vec{\alpha},\vec{\beta}}(\cZ_{a,b},d)$.
\par
The refined Gromov-Witten partition function of $\cX_{a,b}$ is given by
\begin{equation*}
  \begin{aligned}
    Z(\ep_1,\ep_2, Q,x,y) := \exp\Big(\sum_{d=1}^\infty\sum_{g=0}^\infty\sum_{k_1,\dots,k_{a-1}\atop {l_1,\dots,l_{b-1}}}
      N^{\ep_1,\ep_2}_{g,d,\vec{\alpha},\vec{\beta}}Q^d\prod_{i=1}^{a-1}\frac{x_i^{k_i}}{k_i!}\prod_{j=1}^{b-1}\frac{y_j^{l_j}}{l_j!}
    \Big).
  \end{aligned}
\end{equation*}
The refined Gromov-Witten partition function of $\cX_{a,b}$ is given by gluing two refined GW vertex $\tilde{K}_{\mu}(x)$. 
The vertex attached to $p_0$ is given in Theorem \ref{thm:refine-vertex}. On the other side,
the vertex attached to the point $p_\infty$ is obtained by setting
\begin{equation}\label{eqn:infty}
  \begin{aligned}
    &u \rightarrow -u, \quad \tau \rightarrow w, \quad w\rightarrow\tau
    \\
    &\ep_1\rightarrow -\ep_2,\quad \ep_2 \rightarrow -\ep_1, 
  \end{aligned}
\end{equation}
and the framing condition is given by $w = -\ep_2$. (See Figure \ref{fig:opposite})

\begin{figure}[H]

\tikzset{every picture/.style={line width=0.65pt}} %set default line width to 0.75pt        

\begin{tikzpicture}[x=0.60pt,y=0.60pt,yscale=-1,xscale=1]
%uncomment if require: \path (0,424); %set diagram left start at 0, and has height of 424

%Straight Lines [id:da7334901322790686] 
\draw    (200,200) -- (400,200) ;
%Straight Lines [id:da7860048751318542] 
\draw    (400,142) -- (400,252) ;
\draw [shift={(400,254)}, rotate = 270] [color={rgb, 255:red, 0; green, 0; blue, 0 }  ][line width=0.75]    (10.93,-3.29) .. controls (6.95,-1.4) and (3.31,-0.3) .. (0,0) .. controls (3.31,0.3) and (6.95,1.4) .. (10.93,3.29)   ;
\draw [shift={(400,140)}, rotate = 90] [color={rgb, 255:red, 0; green, 0; blue, 0 }  ][line width=0.75]    (10.93,-3.29) .. controls (6.95,-1.4) and (3.31,-0.3) .. (0,0) .. controls (3.31,0.3) and (6.95,1.4) .. (10.93,3.29)   ;
%Straight Lines [id:da14673932587529215] 
\draw    (200,200) -- (200,142) ;
\draw [shift={(200,140)}, rotate = 90] [color={rgb, 255:red, 0; green, 0; blue, 0 }  ][line width=0.75]    (10.93,-3.29) .. controls (6.95,-1.4) and (3.31,-0.3) .. (0,0) .. controls (3.31,0.3) and (6.95,1.4) .. (10.93,3.29)   ;
%Straight Lines [id:da2914220884854658] 
\draw    (200,200) -- (151.28,258.46) ;
\draw [shift={(150,260)}, rotate = 309.81] [color={rgb, 255:red, 0; green, 0; blue, 0 }  ][line width=0.75]    (10.93,-3.29) .. controls (6.95,-1.4) and (3.31,-0.3) .. (0,0) .. controls (3.31,0.3) and (6.95,1.4) .. (10.93,3.29)   ;
%Straight Lines [id:da6546981152430535] 
\draw    (200,200) -- (238,200) ;
\draw [shift={(240,200)}, rotate = 180] [color={rgb, 255:red, 0; green, 0; blue, 0 }  ][line width=0.75]    (10.93,-3.29) .. controls (6.95,-1.4) and (3.31,-0.3) .. (0,0) .. controls (3.31,0.3) and (6.95,1.4) .. (10.93,3.29)   ;
%Straight Lines [id:da7153349339200789] 
\draw    (200,200) -- (157.74,167.99) ;
\draw [shift={(156.14,166.79)}, rotate = 37.14] [color={rgb, 255:red, 0; green, 0; blue, 0 }  ][line width=0.75]    (10.93,-3.29) .. controls (6.95,-1.4) and (3.31,-0.3) .. (0,0) .. controls (3.31,0.3) and (6.95,1.4) .. (10.93,3.29)   ;
%Straight Lines [id:da7926484723398272] 
\draw    (200,200) -- (146.12,191.11) ;
\draw [shift={(144.14,190.79)}, rotate = 9.37] [color={rgb, 255:red, 0; green, 0; blue, 0 }  ][line width=0.75]    (10.93,-3.29) .. controls (6.95,-1.4) and (3.31,-0.3) .. (0,0) .. controls (3.31,0.3) and (6.95,1.4) .. (10.93,3.29)   ;
%Straight Lines [id:da9151267347557622] 
\draw    (400,199.7) -- (357.74,166.99) ;
\draw [shift={(356.14,165.79)}, rotate = 37.14] [color={rgb, 255:red, 0; green, 0; blue, 0 }  ][line width=0.75]    (10.93,-3.29) .. controls (6.95,-1.4) and (3.31,-0.3) .. (0,0) .. controls (3.31,0.3) and (6.95,1.4) .. (10.93,3.29)   ;
%Straight Lines [id:da60735066210283] 
\draw    (400,199.7) -- (346.12,190.11) ;
\draw [shift={(344.14,189.79)}, rotate = 9.37] [color={rgb, 255:red, 0; green, 0; blue, 0 }  ][line width=0.75]    (10.93,-3.29) .. controls (6.95,-1.4) and (3.31,-0.3) .. (0,0) .. controls (3.31,0.3) and (6.95,1.4) .. (10.93,3.29)   ;

% Text Node
\draw (209,138.4) node [anchor=north west][inner sep=0.75pt]    {$-wu$};
% Text Node
\draw (355,139.4) node [anchor=north west][inner sep=0.75pt]    {$-wu$};
% Text Node
\draw (166,246.4) node [anchor=north west][inner sep=0.75pt]    {$(\tau+\frac{1}{b})u$};
% Text Node
\draw (370,237.4) node [anchor=north west][inner sep=0.75pt]    {$\tau u$};
% Text Node
\draw (228,174.4) node [anchor=north west][inner sep=0.75pt]    {$-\frac{u}{b}$};
% Text Node
\draw (118,146.4) node [anchor=north west][inner sep=0.75pt]    {$-\epsilon_{2}u$};
% Text Node
\draw (101,188.4) node [anchor=north west][inner sep=0.75pt]    {$-\epsilon_{1}u$};
% Text Node
\draw (310,149.4) node [anchor=north west][inner sep=0.75pt]    {$-\epsilon_{2}u$};
% Text Node
\draw (302,180.4) node [anchor=north west][inner sep=0.75pt]    {$-\epsilon_{1}u$};

\end{tikzpicture}
\caption{refined vertex attached at $p_\infty$}\label{fig:opposite}
\end{figure}
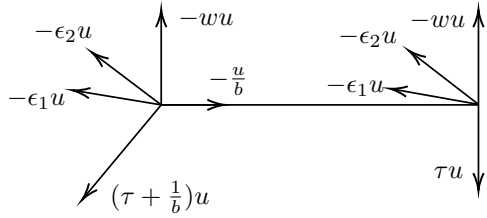

Let $q=e^{\ep_2 u}, s_l=\xi_b^{-1}e^{-\sum_{i=1}^{b-1}\frac{w_b^{-2il}}{b}(w_b^i-w_b^{-i})y_i},l=1,\cdots,b-1$, where $\xi_b = e^{\frac{2\pi i}{b}}, w_b = e^{\frac{\pi i}{b}}$.
Then in the framing condition $w=-\ep_2$, we have the refined GW vertex attached to $p_\infty$:
  \begin{equation}\label{eqn:infty-vertex}
    \begin{aligned}
      \xo{\tilde{K}}(y) \ &= \ G^\bullet_\mu(-\sqrt{-1}\ep_2u;0;y)_b
      \\
      \ &= \ \sqrt{-1}^{|\mu|}(q^{\frac{1}{2}}s_1^{-\frac{1}{b}}\cdots s_{b-1}^{-\frac{b-1}{b}})^{|\mu|}
      \sum_{|\nu|=|\mu|}s_{\nu'}(\tilde{q}_\bullet)\frac{\chi_\nu(\mu)}{z_\mu},
    \end{aligned}
  \end{equation} 
  where $\tilde{q}_{\bullet} = (\tilde{Q},\tilde{Q}s_{b-1},\dots, \tilde{Q}s_1\cdots s_{b-1})$ and $\tilde{Q} = (1,q^{1},q^{2},q^{3},\dots)$. 
  In the smooth case (when $b=1$ and $y = \emptyset$), let $q^{-\rho} = (q^{\frac{1}{2}},q^{\frac{3}{2}},\dots)$, we have
\[
    \xo{\tilde{K}} = \sqrt{-1}^{|\mu|}\sum_{|\nu|=|\mu|}s_{\nu'}(q^{-\rho})\frac{\chi_{\nu}(\mu)}{z_\mu}. 
\]

\begin{remark} 
  We assume $\ep_1>0, \ep_2<0$, and we fix the branch of the square root by defining: $$\sqrt{-\ep_1^2}= \sqrt{-1}\ep_1, \quad \sqrt{-\ep_2^2}=-\sqrt{-1}\ep_2.$$
  This choice ensures consistency when specializing to the unrefined limit $\ep_1=-\ep_2=\ep$.
\end{remark}

\subsection{The gluing rule}
We apply the localization formula in \cite{Liu13} to obtain the gluing rule of refined GW vertex. We introduce the decorated graph $\vec{\Ga}$ in \cite{Liu13}. 
\begin{definition}
Let $n\in \bZ_{\geq 0}$, $\vec{i}= (i_1,\dots, i_n)$, where $-b+1\leq i_j\leq a-1$.
A decorated graph $\vec{\Ga}=(\Ga,\vec{f},\vec{d},\vec{g},\vec{s},\vec{k})$ for $n$-pointed, genus $g$, degree $d$ stable maps to $\cZ_{a,b}$ consists of the following data.
\begin{itemize}
  \item [(1)] $\Ga$ is a compact, connected 1 dimensional CW complex. The set of vertices (resp. edges) in $\Ga$ is denoted by $V(\Ga)$ (resp. $E(\Ga)$). The set of flags of $\Ga$ is defined by
  \[
      F(\Ga) = \{(v,e)\in V(\Ga)\times E(\Ga): v\in e\}.
  \] 
  \item [(2)] The label map $\vec{f}:V(\Ga)\rightarrow \{0,\infty\}$ sends a vertex $v\in V(\Ga)$ to $\{0, \infty\}$. If $v,v'$ are vertices connected by an edge, then $\vec{f}(v)\neq \vec{f}(v')$.
  \item [(3)] The degree map $\vec{d}: E(\Ga)\rightarrow \bZ_{>0}$ sends an edge $e\in E(\Ga)$ to a positive integer $d_e$.
  \item [(4)] The genus map $\vec{g}: V(\Ga)\rightarrow\bZ_{\geq 0}$ sends a vertex $v\in V(\Ga)$ to a nonnegative integer $g_v$.
  \item [(5)] The marking map $\vec{s}:\{1,\dots,n\} \rightarrow V(\Ga)$ is defined if $n>0$.
  \item [(6)] The twisting map $\vec{k}$ sends a flag $(v,e)\in F(\Ga)$ to an element $k_{(v,e)}\in G_v$, a marking $j\in\{1,\dots,n\}$ to an element $k_j\in G_v$, where $G_v = \bZ_a$ (\emph{resp.} $\bZ_b$) if $\vec{f}(v)=0$ (\emph{resp.} $\vec{f}(v)=\infty$). 
\end{itemize}
The above maps satisfy the following constraints:
\begin{itemize}
  \item [(i)] (genus) $\sum_{v\in V(\Ga)}g_v + |E(\Ga)| - |V(\Ga)| + 1 = g$.
  \item [(ii)] (degree) $\sum_{d\in E(\Ga)}d_e = d$.
  \item [(iii)] (compatibility along an edge) Given any edge $e\in E(\Ga)$, let $v,v'\in V(\Ga)$ be its two ends, 
  then $k_{(v,e)}$ and $k_{(v',e)}$ are determined by $d_e$.
  \item [(iv)] (compatibility along a vertex) Given $v\in V(\Ga)$, let $E_v$ be the edges attached at $v$, and $S_v$ be the marking attached at $v$. Then
          \[
              \prod_{e\in E_v}k_{(v,e)}^{-1}\prod_{j\in S_v}k_j = 1.
          \]
  \item [(v)] (compatibility with $\vec{i}$) For each $j=1,\dots,n$, let $v = \vec{s}(j)$, the pair $(p_{\vec{f}(v)},k_j)$ represents a point in the inertia component $\cZ_{i_j}$.
\end{itemize}
\end{definition} 
Let $G_{g, \vec{\alpha}\sqcup -\vec{\beta}}$ be the set of the decorated graphs defined above. According to \cite[Theorem 137]{Liu13}, we have
\begin{equation*}
  \begin{aligned}
    & N_{g,d,\vec{\alpha},\vec{\beta}}^{\ep_1,\ep_2}
    \\
    = & \ \sum_{\vec{\Ga}\in G_{g, \vec{\alpha}\sqcup -\vec{\beta}}(\cZ_{a,b},d)} c_{\vec{\Ga}}\prod_{e\in E(\Ga)}{\bf h}(e)\prod_{(v,e)\in F(\Ga)}{\bf h}(v,e)
    \\
    & \ \cdot \prod_{v\in V(\Ga)}\int_{\overline{\cM}_{g_v,\vec{i}_v}(\cB G_v)}\frac{{\bf h}(v)}{\prod_{e\in E_v}(w_{(v,e)}-\bar{\psi}_{(v,e)}/r_{(v,e)})}. 
  \end{aligned}
\end{equation*}
Here, $G_v = \bZ_a$ or $\bZ_b$. For convenience, we explain the above notations under the assumption $G_v=\bZ_a$. The case of $\bZ_b$ can be obtained by replacing $a\rightarrow b$ and \eqref{eqn:infty}.
\begin{itemize}
  \item [$\bullet$] \emph{Coefficients}. 
  Let $r_{(v,e)} = \frac{|G_v|}{\gcd(d(e),|G_v|)}$, $w_{(v,e)} = \frac{u}{r_{(v,e)}d(e)}$, and 
  $$
    c_{\vec{\Ga}} = \frac{1}{|\Aut(\vec{\Ga})|\prod_{e\in E(\Ga)}d(e)}\prod_{(v,e)\in F(\Ga)}\frac{|G_v|}{r_{(v,e)}}.
  $$
  \item [$\bullet$] \emph{Edge contribution ${\bf h}(e)$}. Let $d=d(e)$, 
  \begin{equation*}
    \begin{aligned}
      {\bf h}(e)
      & = \frac{(\frac{d}{u})^{[\frac{d}{a}]}}{[\frac{d}{a}]!}\frac{(-\frac{d}{u})^{[\frac{d}{b}]}}{[\frac{d}{b}]!}
      \frac{\prod_{j=1}^{[\frac{d}{a}]-\delta_{0,\langle\frac{d}{a}\rangle}}(-(dw+j)u)}{d^{[\frac{d}{a}]-\delta_{0,\langle\frac{d}{a}\rangle}}}
      \\
      & \quad \times  \frac{\prod_{j=1}^{[\frac{d}{b}]-\delta_{0,\langle\frac{d}{b}\rangle}}((d\tau+j)u)}{d^{[\frac{d}{b}]-\delta_{0,\langle\frac{d}{b}\rangle}}}\cdot \frac{1}{\ep_1 u}\frac{1}{\ep_2 u}.
    \end{aligned}
  \end{equation*}
  \item [$\bullet$] \emph{Flag contribution ${\bf h}(v,e)$}.  $${\bf h}(v,e) = (\tau u)(\ep_1 u)(\ep_2 u)(\frac{u}{a})^{\delta_{0,\langle\frac{d(e)}{a}\rangle}}((-w-\frac{1}{a})u)^{\delta_{0,\langle\frac{d(e)}{a}\rangle}}.$$
  \item [$\bullet$] \emph{Vertex contribution ${\bf h}(v)$}.
  \begin{equation*}
    \begin{aligned}
      {\bf h}(v)  = \ &\Lambda^{\vee,U}_{g(v)}(\frac{u}{a})(\frac{u}{a})^{-\delta_v}\cdot \Lambda^{\vee,1}_{g(v)}(\tau u)(\tau u)^{-1}
            \\
            \ \cdot \ &  \Lambda^{\vee,U^{\vee}}_{g(v)}((-w-\frac{1}{a})u)((-w-\frac{1}{a})u)^{-\delta_v}
            \\
            \ \cdot  \ &  \Lambda^{\vee,1}_{g(v)}(-\ep_1 u)(-\ep_1 u)^{-1}\cdot \Lambda^{\vee,1}_{g(v)}(-\ep_2 u)(-\ep_2 u)^{-1}.
    \end{aligned}
  \end{equation*}
\end{itemize}
For any vertex $v$, let $\mu$ denote the profile of the vector $(d(e) : (v,e)\in F(\Ga))$. Define
\begin{equation*}
  \begin{aligned}
    \tilde{I}^0(v) := & \ a^{l(\mu)}\cdot\prod_{i=1}^{l(\mu)}\frac{\prod_{l=1}^{[\frac{\mu_i}{a}]}(\mu_i w+l)}{[\frac{\mu_i}{a}]!(u-\mu_i\bar{\psi}_i)}
    \\
    &\cdot (-1)^{\sum_{(v,e)\in F(\Ga)}[\frac{d(e)}{a}]-\delta_v}\Lambda^{\vee,U}_{g(v)}(\frac{u}{a})\cdot(\frac{u}{a})^{\sum_{(v,e)\in F(\Ga)}\delta_{0,\langle\frac{d(e)}{a}\rangle}-\delta_v}
    \\
    & \cdot \Lambda^{\vee,1}_{g(v)}(\tau u)\cdot (\tau u)^{\val(v)-1}
    \\
    & \cdot \Lambda^{\vee,U^\vee}_{g(v)}((-w-\frac{1}{a})u)\cdot \Big((w+\frac{1}{a})u\Big)^{-\delta_v}
    \\
    & \cdot \Lambda^{\vee,1}_{g(v)}(-\ep_1 u)\cdot (-\ep_1 u)^{\val(v)-1}
    \cdot \Lambda^{\vee,1}_{g(v)}(-\ep_2 u)\cdot (-\ep_2 u)^{\val(v)-1},
    \\
    \tilde{G}^0(v):=  & \ (\ep_1 u)^{-l(\mu)}\int_{\overline{\cM}_{g_v,\vec{i}_v}(\cB\bZ_a)}\tilde{I}^0(v).
  \end{aligned}
\end{equation*}
Similarly, by replacing $a\rightarrow b$ and \eqref{eqn:infty}, we define $\tilde{I}^\infty(v)$ and $\tilde{G}^\infty(v)$ for $v\in V(\Ga)$ with $G_v = \bZ_b$. We have
\begin{equation*}
    N_{g,d,\vec{\alpha},\vec{\beta}}^{\ep_1,\ep_2} = \sum_{\vec{\Ga}\in G_{g,\vec{\alpha}\sqcup-\vec{\beta}}(\cZ_{a,b},d)}\frac{1}{|\Aut(\vec{\Ga})|}\prod_{v\in V(\Ga)}\tilde{G}^{\vec{f}(v)}(v)\prod_{e\in E(\Ga)}d_e.
\end{equation*}
Under the framing condition, $\tau = \ep_1, w=-\ep_2$, we have 
\begin{align*}
      Z(\ep_1,\ep_2,Q,x,y) &= \sum_{\mu}G^\bullet_{\mu}(\sqrt{-1}\ep_1 u;w;x)_a(-1)^{l(\mu)-|\mu|}Q^{|\mu|}z_\mu G^\bullet_\mu(-\sqrt{-1}\ep_2 u;\tau;y)_b
    \\
      &= \sum_{\mu}G^\bullet_{\mu}(\sqrt{-1}\ep_1 u;0;x)_a(-1)^{l(\mu)-|\mu|}Q^{|\mu|}z_\mu G^\bullet_\mu(-\sqrt{-1}\ep_2 u;0;y)_b\\
      &=\prod_{i,j\geq 0}\prod_{1\leq k\leq a\atop{1\leq l\leq b}}
      \Big(1 - Qt^{i+\frac{1}{2}}q^{j+\frac{1}{2}}\cdot \hat{q}_k\hat{s}_l\cdot q_1^{-\frac{1}{a}}\cdots q_{a-1}^{-\frac{a-1}{a}}
      \cdot s_1^{-\frac{1}{b}}\cdots s_{b-1}^{-\frac{b-1}{b}}\Big)
\end{align*}
where the second equality follows from Theorem \ref{key} and \eqref{eqn:Comp} \eqref{eqn:Init} \eqref{eqn:minus-phi}, the third 
equality follows from Theorem \ref{thm:refine-vertex} and \eqref{eqn:infty-vertex}, and 
\begin{align*}
	\hat{q}_k &= \left\{
	\begin{aligned}
		&q_k\cdots q_{a-1}, && 1\leq k\leq a-1,
		\\
		&1, && k=a.
	\end{aligned}\right.
	,
	\ \text{and}\
	\hat{s}_l = \left\{
	\begin{aligned}
		&s_l\cdots s_{b-1}, && 1\leq l\leq b-1,
		\\
		&1, && l=b.
	\end{aligned}\right.
\end{align*}

\subsection{The resolved conifold}
When $\cX$ is the resolved conifold, i.e. $a=b=1$, we have
\begin{equation*}
  \begin{aligned}
      Z(\ep_1,\ep_2, Q) &= \sum_{\mu} Q^{|\mu|}z_\mu \sum_{|\nu|=|\eta|=|\mu|} s_{\nu'}(t^{-\rho})s_{\eta'}(q^{-\rho})\frac{\chi_\nu(\mu)}{\sqrt{-1}^{l(\mu)}z_\mu}\frac{\chi_\eta(\mu)}{\sqrt{-1}^{l(\mu)}z_\mu}
      \\
      &= \sum_\nu (-Q)^{|\nu|}s_{\nu}(t^{-\rho})s_{\nu'}(q^{-\rho})
      \\
      &= \prod_{i,j=0}^\infty \Big(1 - Qt^{i+\frac{1}{2}}q^{j+\frac{1}{2}}\Big)
      \\
      &= \exp\left\{-\sum_{d=1}^\infty \frac{Q^d}{d(t^{\frac{d}{2}}-t^{-\frac{d}{2}})(q^{\frac{d}{2}}-q^{-\frac{d}{2}})}\right\}.
  \end{aligned}
\end{equation*}

\nocite{*}
\bibliographystyle{plain}
\bibliography{refs}

\end{document}